\documentclass[final,12pt,3p,times]{elsarticle}

\usepackage{etoolbox}
\patchcmd{\MaketitleBox}{\footnotesize\itshape\elsaddress\par\vskip36pt}{\footnotesize\itshape\elsaddress\par\parbox[b][36pt]{\linewidth}{\vfill\hfill\textnormal{February 6, 2024}\hfill\null\vfill}}{}{}%

\usepackage{import}
\usepackage{thesis}

\usepackage{amssymb}

\usepackage{lineno}

\begin{document}


\begin{frontmatter}



\title{Corotational modeling and NURBS-based kinematic constraint implementation in three-dimensional vehicle-track-structure interaction analysis}


\author[inst1]{Maria Fedorova\fnref{fnote1}} 
\ead{mariafed@buffalo.edu}
\author[inst1]{M.V. Sivaselvan\corref{cor1}} 
\ead{mvs@buffalo.edu}
\cortext[cor1]{Corresponding author}
\fntext[fnote1]{Current address: LARSA, Inc., 68 S Service Rd, Ste 100, Melville, 11747, NY, USA}

\affiliation[inst1]{organization={Department of Civil, Structural, and Environmental Engineering, University at Buffalo},
            addressline={Ketter Hall}, 
            city={Buffalo},
            postcode={14260}, 
            state={NY},
            country={USA}}

\begin{abstract}
An algorithm for three-dimensional dynamic vehicle-track-structure interaction (VTSI) analysis is described in this paper. The algorithm is described in terms of bridges and high-speed trains, but more generally applies to multibody systems coupled to deformable structures by time-varying kinematic constraints. Coupling is accomplished by a kinematic constraint/Lagrange multiplier approach, resulting in a system of index-3 Differential Algebraic Equations (DAE). Three main new concepts are developed. (i) A corotational approach is used to represent the vehicle (train) dynamics. Reference coordinate frames are fitted to the undeformed geometry of the bridge. While the displacements of the train can be large, deformations are taken to be small within these frames, resulting in linear (time-varying) rather than nonlinear dynamics. (ii)~ If conventional finite elements are used to discretize the track, the curvature is discontinuous across elements (and possibly rotation, too, for curved tracks). This results in spurious numerical oscillations in computed contact forces and accelerations, quantities of key interest in VTSI. A NURBS-based discretization is employed for the track to mitigate such oscillations. (iii) The higher order continuity due to using NURBS allows for alternative techniques for solving the VTSI system. First, enforcing constraints at the acceleration level reduces an index-3 DAE to an index-1 system that can be solved without numerical dissipation. Second, a constraint projection method is proposed to solve an index-3 DAE system without numerical dissipation by correcting wheel velocities and accelerations. Moreover, the modularity of the presented algorithm, resulting from a kinematic constraint/Lagrange multiplier formulation, enables ready integration of this VTSI approach in existing structural analysis and finite element software. A simplified vehicle-bridge model is first introduced that exemplifies numerical challenges. This model is used throughout the paper to motivate various algorithmic choices. Finally, a numerical example consisting of realistic vehicle and bridge models is presented. 
\end{abstract}




\end{frontmatter}


\section{Introduction} \label{sec:intro}

Traditional static analysis with moving loads is not sufficient for many bridges supporting high-speed trains. Typically used codes of practice \cite{CaliforniaHSR,Eurocode1991_2} require dynamic analysis for such structures. It can be performed either with a series of moving loads representing the dynamic effects of multiple trains and axles or by solving a coupled system of vehicle, track, and structure equations of motion. The latter approach constitutes vehicle-track-structure interaction (VTSI) analysis. According to the California High-Speed Train Project documentation \cite{CaliforniaHSR}, VTSI analysis is required for some structures that do not meet the requirements of a simplified analysis. \cite{FedorovaMag2019}

There are numerous ways to model dynamic interaction between trains and track-structure subsystems (see, for example, References \cite{Yang2004,Shabana2007,Zhai2019} for an overview of available approaches).
As was discussed by Fedorova and Sivaselvan \cite{Fedorova2017}, one of the drawbacks of many existing algorithms is the way the governing equations of motion are solved. The equations for the train and bridge subsystems are often combined together into a single equation to obtain the solution \cite{neves2012direct,Xia2005,Dimitrakopoulos2015,Salcher2015,Arvidsson2019,Bettinelli2023,Glatz2021,Greco2018,Mosayebi2021,Gou2018,Xiao2020,Xu2019,Xiao2019,Eroglu2022,Zhou2023_VSD,Yu2018,He2020,Wang2019,Li2022,Zhang2023,Shi2022,Chen2021_VSD,Chen2022,Montenegro2020,Xin2020,Su2022}, or the two sets of equations are solved separately employing iterative procedures \cite{Diana1989,Kwark2004,Zhang2013,TiconaMelo2018,TiconaMelo2020,Xu2022_VSD,Lu2020,Liu2022,Li2021,Zhou2023,Stefanidou2022}. While the former approach makes it difficult to integrate with existing structural analysis software, the latter requires significant computational effort. Another possible approach is to solve the system of equations using a complementary kinematic constraint condition \cite{sivaselvan2014,Zhu2015,Fedorova2017}. This approach facilitates integration into existing software with minimum intervention, and, moreover, modeling of contact separation between the wheels and the bridge becomes more straightforward. Furthermore, wheel-rail contact forces are direct outputs of the computation. 
 
This paper builds on an algorithm for two-dimensional vehicle-track-structure interaction (VTSI) analysis \cite{Fedorova2017}. There are three main new features of the approach presented in this paper: 1) a corotational approach, wherein train dynamics is described with reference to coordinate frames fitted to the bridge geometry; deformations within these frames can be taken to be small, resulting in linear (time-varying) dynamics; 2) NURBS-based discretization of the vehicle path to provide sufficient inter-element continuity, and prevent spurious oscillations in contact forces and accelerations, quantities of principal interest in VTSI; 3) two alternative ways of solving index-3 system of differential-algebraic equations (DAE) that are possible due to higher continuity provided by NURBS: enforcement of constraints at acceleration level, which reduced the system to a index-1 DAE, and constraint projection method, which allows for solving an index-3 system without numerical dissipation by correcting wheel velocities and accelerations. Moreover, while the algorithm is applied specifically to VTSI analysis in the current paper, it can generally be used in multibody systems coupled to deformable structures by time-varying kinematic constraints. 

This paper is organized as follows. Section \ref{sec:corotational} introduces the corotational approach used to describe vehicle kinematics. A numerical example is then introduced (Section \ref{sec:numExample}) and used as a running thread throughout the paper to illustrate the essential features of the algorithm and describe numerical issues and resolutions. These numerical issues are described in Section \ref{sec:problemStatement}, associated with solving VTSI equations as an index-3 DAE and using finite elements to model the track. To resolve the numerical issues, NURBS-based discretization of the track is employed in Section \ref{sec:nurbs}. In Section \ref{sec:algorithm}, we describe three ways of solving VTSI equations: (i) first, a Generalized-$\alpha$ scheme is used to obtain desired numerical damping and mitigate the spurious oscillations; (ii) second, the advantage of using NURBS for modeling the track is demonstrated by explicitly enforcing constraints at acceleration level and avoiding the use of numerical dissipation; (iii) third, we discuss a simple constraint projection method that allows solving VTSI as an index-3 DAE system without numerical dissipation by non-iterative correction of obtained velocities and accelerations at each time step. The advantages and disadvantages of the proposed methods in terms of practical implementation are discussed. While some of the methods are easy to implement, they may not be well suited for existing structural analysis software. Further numerical results for a simplified model, as well as sample results for a realistic model, are included in Section \ref{sec:moreNumResults}. The key findings are then summarized in Section~\ref{sec:conclusions}.

\section{Corotational approach} \label{sec:corotational}

\begin{figure}[h]
\centering
\captionsetup{justification=centering}
\includegraphics[scale=0.8]{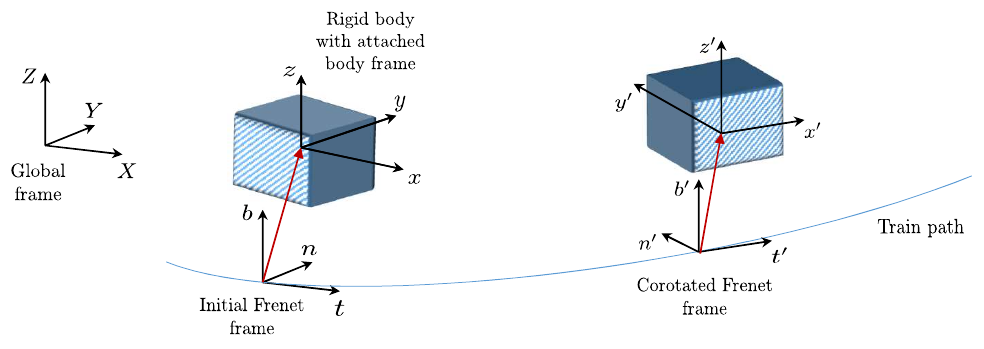}
\caption{The corotational approach in VTSI: motion of a body is described with respect to a moving Frenet frame. Displacements and rotations of the Frenet frame can be large, but displacements and rotations relative to the Frenet frame are small. Here, $XYZ$ is a Global frame, $tnb$ is a Frenet frame, and $xyz$ is a rigid body frame. The apostrophe indicates the frames moved to a new, corotated position.}
\label{fig:frames}
\end{figure}

The corrotational approach comes from continuum mechanics and the finite element method (see, for example, Felippa et al. \cite{Felippa2005} for an overview of the approach and \cite{You2024} for an example of practical implementation). The key idea of this approach is to allow the displacements and rotations of a reference coordinate frame to be large while keeping deformations within the frame small. In the VTSI context, cars of trains transversing curved paths undergo large rigid body motion, but deformations, for instance, of the suspension system, are small. It is, therefore, convenient to express the car dynamics in local corotated coordinate frames fit to the curved profile of the path. In this paper, reference Frenet frames are fit to the path of the vehicle (Figure \ref{fig:frames}). 
At the start of the analysis, the vehicle's initial position at the beginning of the path is described with respect to the initial Frenet frame (see Figure \ref{fig:frames}). 
Later, as the vehicle moves, the deformations within each car of the train are computed in a Frenet frame corresponding to the position of the car body's center of gravity along the curved path of the bridge. Since the train longitudinal dynamics are not modeled in this study, longitudinal connections between the cars (couplers) are not modeled either. The displacements and rotations of the car parts inside the Frenet frame are taken to be small. Each part is taken to be a rigid body, forming a system of rigid bodies (see Figure \ref{fig:modelSimple} for a simplified vehicle model and Figure \ref{fig:trainModelReal} for a realistic train model). As opposed to trajectory coordinates often used in railroad vehicle formulations (for example, Shabana et al. \cite{Shabana2007}), the corotational approach requires deformations within the reference frames to be small, resulting in linear dynamics, which is appropriate in practical situations, but simplifies the VTSI formulation. 

The path taken by a train is a curve in 3D.
Let $\xhat:[0,L] \rightarrow \mathbb{R}^3$ be such a curve parametrized by arc length~$s$ (in Section \ref{sec:nurbs}, this curve will be represented using NURBS). Using parameter $\xi = \hat{\xi}(s)$, this curve can be reparameterized as $\xtilde:[\xi_a,\xi_b] \rightarrow \mathbb{R}^3$, where the end points $\xi_a$ and $\xi_b$ correspond to the start and end of the path. A local Frenet frame $\{\tF,\nF,\bF \}$ can then be defined at each point $\xtilde(\xi)$ of the curve as \cite{Farin1997,tapp2016differential}

\begin{equation}
    \textbf{t} = \frac{\xtilde'}{\left\| \xtilde' \right\|},\
    \textbf{n} = \textbf{b} \times \textbf{t},\ 
    \textbf{b} = \frac{\xtilde' \times \xtilde''}{\left\| \xtilde' \times \xtilde'' \right\|}
\end{equation}

\noindent where a prime denotes a derivative with respect to parameter $\xi$; $\textbf{t}, \textbf{n},$ and $ \textbf{b}$ are the tangent, normal, and binormal vectors. 
The Frenet frame is attached to the curve and serves both as a reference frame to describe the train dynamics (Section \ref{sec:trainSimple}) and as a basis for the geometry and solution space of equations of motion (Section \ref{sec:nurbs}) of a beam representing a simplified bridge/track (Section \ref{sec:bridgeSimple}) or rails in the realistic model (Section \ref{sec:bridgeReal}).

Having introduced the corotational approach in the VTSI context, we next develop the numerical examples that are used to motivate the algorithmic choices.

\section{Numerical example to motivate algorithmic choices} \label{sec:numExample}

\subsection{Simplified vehicle model} \label{sec:trainSimple}

In developing our VTSI algorithm, a number of choices are encountered along the way, such as the use of numerical damping to mitigate the spurious oscillations (Section \ref{sec:genAlpha}), explicit enforcement of the acceleration constraints (Section \ref{sec:accConstr}), and the employment of a simple constraint projection method (Section \ref{sec:constrCorrestion}). These choices are best explained by means of an example. Here, a simplified model is introduced for this purpose. The model is used to illustrate essential features of the three-dimensional VTSI algorithm and numerical issues associated with solving an index-3 DAE system. The model consists of two rigid bodies representing a wheel and a car and a spring and dashpot representing a vehicle suspension, as shown in Figure \ref{fig:modelSimple}. The bridge itself is represented in the simplified model by one curved beam, without explicit rails, sleepers, ballast etc. (Section \ref{sec:bridgeSimple}). Application of the algorithm to realistic models including such features is demonstrated in Section \ref{sec:moreNumResults:real}.

\subsubsection{Kinematics}

\begin{figure}[h]
\centering
\captionsetup{justification=centering}
\includegraphics[scale=0.9]{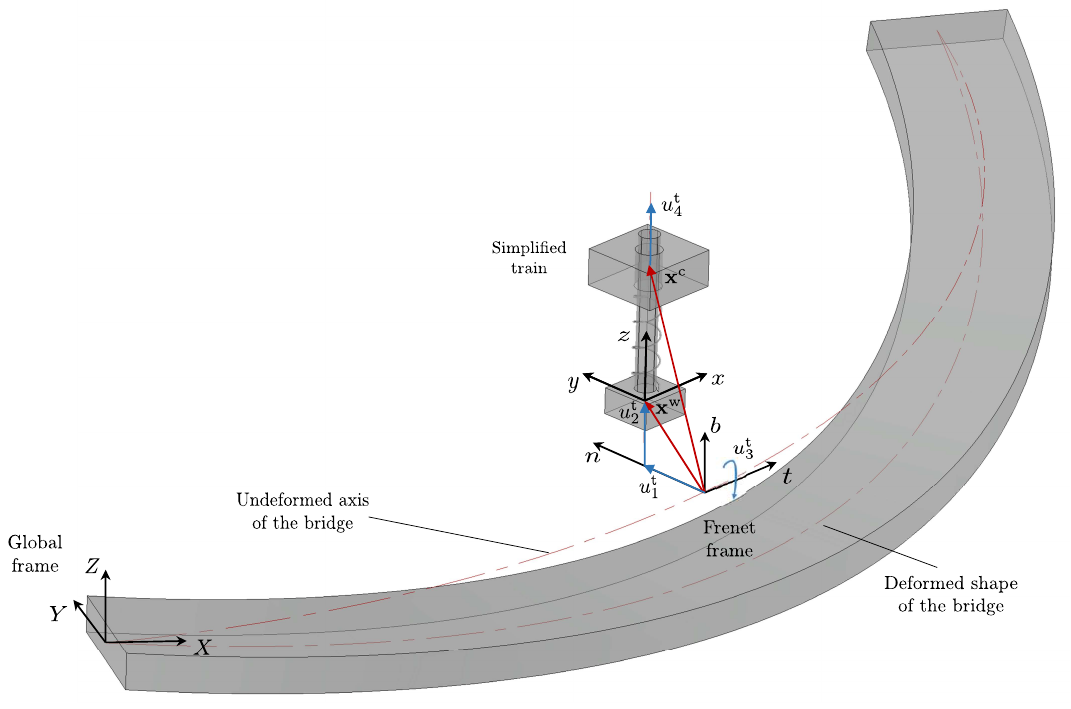}
\caption{Bridge curved in plan and a simplified vehicle model (generalized coordinates of the vehicle model with respect to the Frenet frame are shown in \textcolor{blue}{blue}; position vectors of the vehicle parts are shown in \textcolor{red}{red})}
\label{fig:modelSimple}
\end{figure}

The simplified vehicle model has 4 degrees of freedom. Therefore, four generalized coordinates are required to describe the motion of the model: the transverse $\ut_1$ and vertical $\ut_2$ displacements of the wheel, rolling of the vehicle $\ut_3$, and vertical displacement of the car $\ut_4$. The superscript "t" stands for "train". Figure \ref{fig:modelSimple} shows the simplified vehicle in the Frenet frame attached to the \textit{undeformed} central axis of the bridge.

The positions of the simplified wheel (\xw) and car (\xc) relative to the global frame are then
 
\begin{align} \label{eq:t:pos}
\xw &= \xF + \RF \left( {\begin{array}{*{20}{c}}
    0\\
    {\ut_1}\\
    {\ut_2}
    \end{array}} \right) = 
\xF + \RF {\begin{bmatrix}{}
    0 & 0 & 0 & 0\\
    1 & 0 & 0 & 0\\
    0 & 1 & 0 & 0
    \end{bmatrix}} 
\left( {\begin{array}{*{20}{c}}
    \ut_1\\
    \ut_2\\
    \ut_3\\
    \ut_4
    \end{array}} \right) = 
\xF + \RF \Tw \utB \nonumber
\\
\xc &= \xF + \RF \left( {\begin{array}{c}
    0\\
    {\ut_1 - \ut_3 {l_0}}\\
    {\ut_4 + {l_0}}
    \end{array}} \right) = 
\xF + \RF {\begin{bmatrix}{}
    0 & 0 & 0 & 0\\
    1 & 0 & -l_0 & 0\\
    0 & 0 & 0 & 1
    \end{bmatrix}} 
\left( {\begin{array}{{c}}
    \ut_1\\
    \ut_2\\
    \ut_3\\
    \ut_4
    \end{array}} \right) + \RF 
\left( {\begin{array}{{c}}
    0\\
    0\\
    l_0
    \end{array}} \right) \\ &=
\xF + \RF \Tc \utB + \RF \left( {\begin{array}{c}
    0\\
    0\\
    l_0 
    \end{array}} \right) \nonumber
\end{align}
where ${ \left( \xF, \RF \right) }$ is the specification of the configuration of the Frenet frame relative to the global frame: \xF\ is the position vector of the origin of the frame and \RF\ is an orientation of the frame.

The velocities of the wheel and the car in the global frame are, differentiating \eqref{eq:t:pos},

\begin{equation}  \label{eq:t:vel}
\begin{aligned}
\xwdot &= \xFdot + \RF  \omegaFhat  \Tw \utB 
+ \RF \Tw \utdotB\\
\xcdot &= \xFdot + \RF  \omegaFhat \Tc \utB 
+ \RF \Tc \utdotB 
+ \RF  \omegaFhat 
    \left( {\begin{array}{c}
    0\\
    0\\
    l_0
    \end{array}} \right)
\end{aligned}
\end{equation}
where $\omegaFhat = {\RF}^T \RFdot$ is the matrix representing the angular velocity of the Frenet frame \cite{murray1994}. The matrix \omegaFhat\ is skew-symmetric by definition and serves to map $\mathbb{R}^3 \mapsto \mathbb{R}^{3 \times 3}$, that is, for the angular velocity vector $\omegaFB = {\left( \omegaF_1\ \omegaF_2\ \omegaF_3 \right) }^T$, the matrix \omegaFhat is defined as

\begin{equation}
\omegaFhat = 
\left[ {\begin{array}{*{20}{c}}
0 & -\omegaF_3 & \omegaF_2\\
\omegaF_3 & 0 & -\omegaF_2\\ 
-\omegaF_2 & \omegaF_1 & 0 \end{array}} \right]
\end{equation}
Using properties of cross product operation between vectors \cite{murray1994}, we then can write

\begin{equation}
\omegaFB \times {\bf{x}} = \omegaFhat {\bf{x}}
\end{equation}

Having described vehicle kinematics, we proceed to deriving the equations of motion.

\subsubsection{Equations of motion}

The motion of the system is described by employing Lagrangian mechanics and using the Euler-Lagrange equations of motion for a constrained system \cite{murray1994} (summation over the index $j$ is assumed):

\begin{equation} \label{eq:euler}
\frac{d}{{dt}}\left( {\frac{{\partial L}}{{\partial {{\dot u}_i}}}} \right) - \frac{{\partial L}}{{\partial {u_i}}} + \frac{\partial C_j}{\partial u_i} \lambda_j = 0
\end{equation}
where ${L}$ is a Lagrangian for an unconstrained system, ${u_i}$ are generalized coordinates of the model, $C_j(u)=0$ is a constraint \eqref{eq:eqofmotion:constr}, and $\lambda_j$ is the corresponding Lagrange multiplier. The Lagrangian is constructed from the kinetic energy, ${T}$, and the potential energy, ${V}$, as follows

\begin{equation} \label{eq:lagrangian}
L = T - V
\end{equation}

By substituting \eqref{eq:lagrangian} into \eqref{eq:euler}, and taking into account the fact that the potential energy ${V}$ is not a function of time or velocity, we obtain

\begin{equation} \label{eq:euler2}
\frac{d}{{dt}}\left( {\frac{{\partial T}}{{\partial {{\dot u}_i}}}} \right) - \frac{{\partial T}}{{\partial {u_i}}} + \frac{{\partial V}}{{\partial {u_i}}} + \frac{\partial C_j}{\partial u_i} \lambda_j= 0\\
\end{equation}
The kinetic and potential energies of the simplified train model are given by

\begin{equation} \label{eq:t:energy}
\begin{gathered}
T = \frac{1}{2} \mw \xwdotT \xwdot + \frac{1}{2} \mc \xcdotT \xcdot + \frac{1}{2}\left( \Iw + \Ic \right) {\dot u_3}^{{\text{t}^2}}\\
V = \frac{1}{2} \ks (\ut_4 - \ut_2)^2 + \mw g \left( x^{\text {w}}_3 - x^{\text {w}}_3 \left( 0 \right) \right) + \mc g\left( x^{\text {c}}_3 - x^{\text {c}}_0 \left( 0 \right) \right)\\
\end{gathered}
\end{equation}
where \mw, \mc, \Iw, and \Ic\ are masses and moments of inertia of the simplified wheel and car, and \ks\ is the stiffness of the suspension spring. 

By substituting \eqref{eq:t:vel} into \eqref{eq:t:energy}, and from the Euler-Lagrange equation \eqref{eq:euler2}, we obtain the equations of motion of the model

\begin{equation} \label{eq:t:eqofmotion}
\Mt \utddotB + \Ct(t) \utdotB + \Kt(t) \utB + {\Ltr}^\T \lambdaB(t) = \Pt(t)
\end{equation}
where the train matrices are stated explicitly in \ref{sec:appendix:train}. The complete derivation of Equation \eqref{eq:t:eqofmotion} can be found in Fedorova \cite{FedorovaThesis2017}.

In the numerical example used throughout this paper to motivate algorithmic choices, the vehicle properties are partially adapted from a numerical example by Yang et al. \cite{Yang2004}. The element representing a car of the simplified vehicle model has the mass ${\mc=41750\ \si{kg}}$, mass moment of inertia $\Ic_x = 23.2 \times 10^3\ \si{kg \cdot m^2}$, and its center of gravity located at $1.37\ \si{m}$ above the wheel. 
The wheel element has a mass equal to the mass of four wheelsets ${\mw = 1780 \times 4 = 7120\ \si{kg}}$, and mass moment of inertia $\Ic_w = 1.14 \times 10^3\ \si{kg \cdot m^2}$. The listed mass moments of inertia  $\Ic_x$ and $\Iw_x$ belong to the Shinkansen (SKS) series 300 vehicle model\cite{Yang2004}. The stiffness of the suspension spring is the equivalent stiffness of all primary and secondary suspension springs of a train car
${\ks = 2 / \left( \frac{1}{4 k_\text{prim}} + \frac{1}{2 k_\text{sec}}\ \right)  = 865.6\ \si{kN/m}}$. Here, two bogies per car and two wheelsets per bogie are assumed, with two springs per each suspension. The initial length of the spring is ${l_0 = 1.37\ \si{m}}$ and corresponds to the position of the car body's center of gravity. The speed of the train ${v^\text{t}}$ is constant and equals to ${360\ \si{km/h} = 100\ \si{m/s}}$, which equals the maximum speed of the real trains used in Eurocode EN 1991-2:2003 \cite{Eurocode1991_2}.

\subsection{Simplified bridge model} \label{sec:bridgeSimple}

\begin{figure}[h] 
\centering
\setcounter{subfigure}{0}
\captionsetup[subfigure]{justification=centering}
\begin{tabular}{l@{\hskip 0.8cm}r}
\subfloat[Simplified bridge model (Section~\ref{sec:bridgeSimple})] 
    {\label{fig:modelPlanSimple} 
    \includegraphics{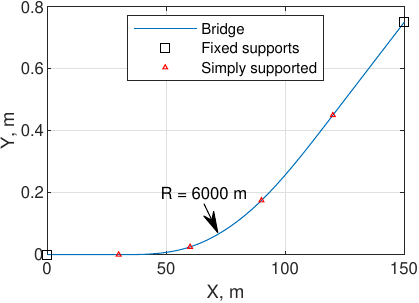}} &
\subfloat[Realistic bridge/track model, ballast elements are not shown     (Section~\ref{sec:bridgeReal})] 
    {\label{fig:modelPlanReal} \includegraphics{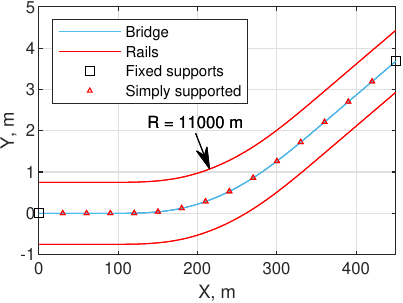}}
\end{tabular}
\caption{Paths in simplified and realistic bridge models used in numerical examples}
\label{fig:modelPlan}
\end{figure}

A simplified bridge model is first used for clarity of presentation; a realistic model is used in the numerical example in Section \ref{sec:moreNumResults:real}. The simplified model does not distinguish between the track and the deck (Figure \ref{fig:modelSimple}). The bridge can be modeled using standard FEM elements or using NURBS (Section \ref{sec:nurbs}). The numerical example used throughout this paper employs a simplified bridge model that consists of two straight spans, two transition spans, and one curved span in the middle of the structure (Figure \ref{fig:modelPlanSimple}). The bridge is fixed at the ends and simply-supported in the middle supports. 

The material and geometric properties of the bridge are partially adapted from a straight bridge model by Yang et al. \cite{Yang2004}. The length of each span is ${30\ \si{m}}$. According to the Technical Memorandum TM 2.1.2 \cite{CaliforniaHSR_Track_TM_2_1_2} by the California High-Speed Train Project, the "exceptional" (most minimal) track radius for speed ${355\ \si{km/h}}$ should be ${6000\ \si{m}}$. At the same time, the "desirable" radius for such speed is ${10700\ \si{m}}$. Considering the train speed of ${360\ \si{km/h}}$ (Section \ref{sec:trainSimple}), for the simplified bridge model the radius of ${6000\ \si{m}}$ is chosen. Transition spans are modeled between the straight and curved spans. The radius of the transition spans is ${0\ \si{m}}$ at the beginning of the span and $6000\ \si{m}$ at the end, with linear variation along the length of the span. Track with the radius of ${11000\ \si{m}}$ is used in the realistic bridge model (Figure \ref{fig:modelPlanReal}, Sections \ref{sec:bridgeReal} and \ref{sec:moreNumResults:real}).

The following material properties are used: Young's modulus ${E^b=28.25\ \si{GPa}}$, modulus of rigidity ${G^b=10^3\ \si{GPa}}$, torsional constant ${J = 15.65\ \si{m^4}}$, cross-section area ${A = 7.73\ \si{m^2}}$, cross-section second moments of inertia ${I_y^b=7.84\ \si{m^4}}$ and ${I_z^b=74.42\ \si{m^4}}$, and the mass per unit length (including mass of the ballast) ${41740\ \si{kg/m}}$. The standard finite element (constructed with Hermitian cubic shape functions) or NURBS of degrees $p=3$ and $p = 5$ are used. 

\subsection{Realistic vehicle model} \label{sec:trainReal}

While the simplified vehicle model (Figure \ref{fig:modelSimple}) in Section \ref{sec:trainSimple} is meant to illustrate the essential features of the proposed algorithm, a realistic train model is also presented for practical three-dimensional VTSI analysis. In this paper, the use of the realistic model is restricted to Section \ref{sec:moreNumResults:real}. Additional analysis results and the model derivations are available in Fedorova \cite{FedorovaThesis2023}.

The vehicle equations of motion describe the behavior of rigid bodies (representing train cars, bogies, and wheelsets) and include damping and stiffness terms corresponding to the suspension system that connects these rigid bodies 

\begin{equation} \label{eq:rb:assembly} 
\begin{gathered}
    \assembly \limits_{i = 1}^{n} \left[ \Mrbi \qddot^i + \Crbi  (t) \qdot^i + \Krbi  (t) {\qB}^{i} \right] + 
    \assembly \limits_{j = 1}^{l} 
        \left[
        \Csuspj 
            \left(
            \begin{array}{@{}c@{}}
                {\bf \dot q}^{j_1}\\ 
                {\bf \dot q}^{j_2}
            \end{array} 
            \right)
        \right]
    +\\
    \assembly \limits_{k = 1}^{m} 
        \left[
        \Ksuspk 
            \left(
            \begin{array}{@{}c@{}}
                {\bf \dot q}^{k_1}\\ 
                {\bf \dot q}^{k_2}
            \end{array} 
            \right)
        \right]
    +
    \assembly \limits_{p = 1}^{q} 
        \left[
        \Lrwp 
        \lambdaB^{\text{rw},\ p} (t)
        \right]   
    + 
    \assembly \limits_{r = 1}^{s} 
        \left[
        \Ltrr 
        \lambdaB^{\text{t},\ r} (t)
        \right]
    =
    \assembly \limits_{i = 1}^{n} \left[ \Prbi  (t) \right]
\end{gathered}
\end{equation}
where $\Mrbi$, $\Crbi$, $\Krbi$ are the mass, damping, and stiffness matrices, $\Prbi$ is the load vector, and ${\qB}^{i}$ is the displacement vector of $i$-th rigid body; $\Csuspj$ is a damping matrix and ${\bf \dot q}^{j_1}$ is a velocity vector of the first joint of $j$-th suspension dashpot; $\Ksuspk$ is a stiffness matrix and ${\bf \dot q}^{k_1}$ is a velocity vector of the first joint of $k$-th suspension spring; $\Lrwp$ and $\lambdaB^{\text{rw},\ p}$ are wheel-rail constraint matrix and contact forces of $p$-th wheelset; $\Ltrr$ and $\lambdaB^{\text{t},\ r}$ are kinematic constraint matrix and corresponding constraint forces of $r$-th rigid body; $n$ is the number of rigid bodies (train cars, bogies, and wheelsets), $l$ is the number of dashpots, $m$ is the number of springs, $q$ is the number of wheelsets, $s$ is the number of rigid bodies subjected to kinematic constraints; $\assembly$ is the assembly operator. 

\begin{figure}[h!]
\centering
\captionsetup{justification=centering}
\includegraphics[scale=0.8]{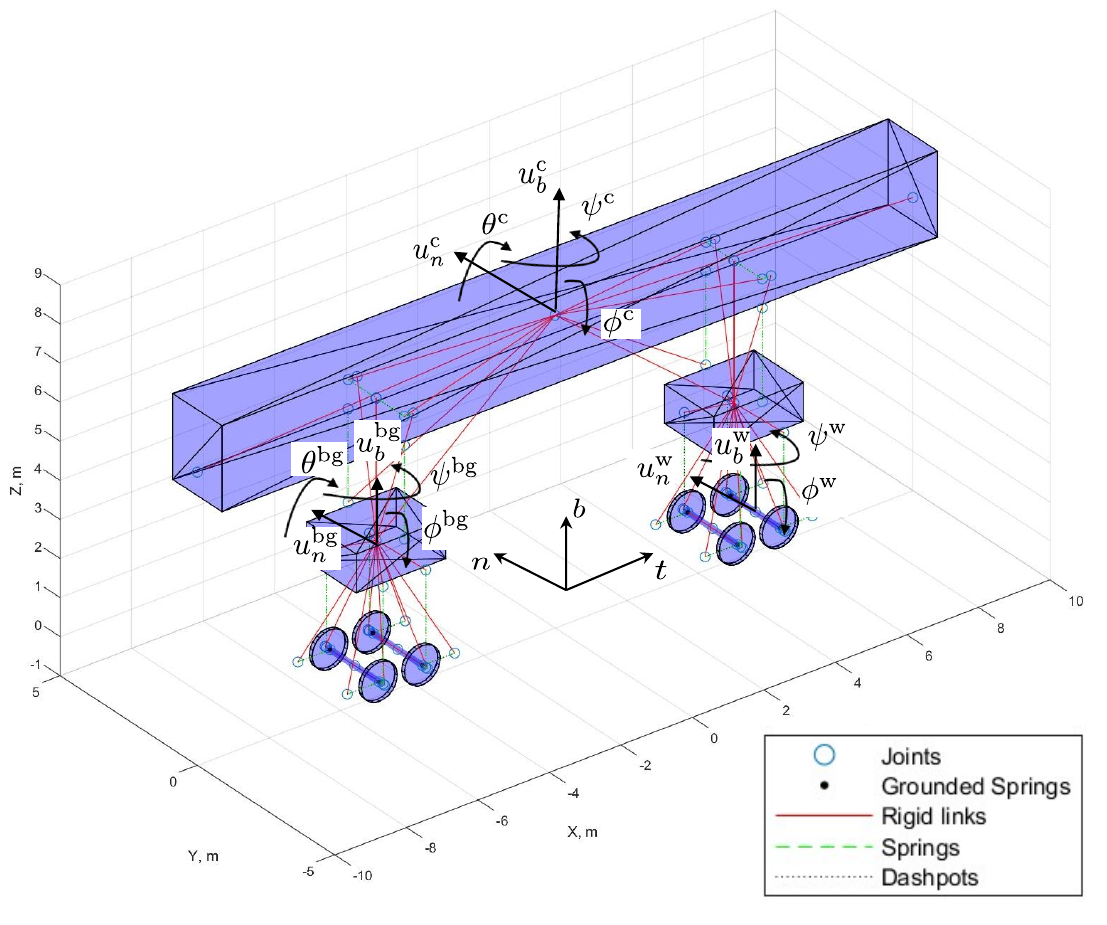}
\caption{An example of a realistic train model (one car is shown)}
\label{fig:trainModelReal}
\end{figure}

Equation \ref{eq:rb:assembly} can be written in a more concise form as
\begin{equation} \label{eq:trainReal} 
\begin{gathered}
    \Mt \qtddot + \Ct (t) \qtdot+ \Kt  (t) \qt + \Lrw \lambdaB^{\text{rw}} (t) + \Ltr \lambdaB^{\text{t}} (t) = \Pt (t) 
\end{gathered}
\end{equation}
where \Lrw\ is a kinematic constraint matrix used to constrain wheels to the rails and $\lambdaB^{\text{rw}}$ is a vector of corresponding contact forces; \Ltr\ is kinematic constraint matrix used to restrict longitudinal motions of rigid bodies and rotations of wheelsets with respect to transverse axis, $n$, and  $\lambdaB^{\text{t}}$ is a vector of corresponding constraint forces. An example train car with two bogies and two wheelsets per bogie then has a total of 31 degrees of freedom, as shown in Figure \ref{fig:trainModelReal}.

\subsection{Realistic bridge model} \label{sec:bridgeReal}

\begin{figure}[h]
\centering
\captionsetup{justification=centering}
\includegraphics[scale=1]{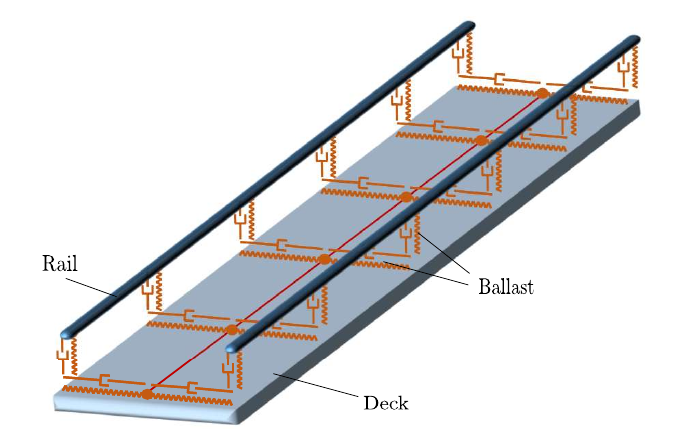}
\caption{Schematic of a realistic bridge and track model}
\label{fig:bridgeReal}
\end{figure}

A realistic bridge model represents the deck and rails as distinct components and accounts for the ballast (see Figure \ref{fig:modelPlanReal}, where ballast elements are not shown, and Figure \ref{fig:bridgeReal}). When the analysis is implemented in finite element analysis software, the deck is modeled with standard finite elements, and the ballast is modeled with springs acting in vertical and transverse directions. The rails can be modeled either with NURBS elements (Section \ref{sec:nurbs}) or with standard finite elements. The rationale behind modeling the rails with NURBS is given in Section \ref{sec:numericalChallenges}. The Kinematic constraints are used to connect the rails to the ballast. This approach provides the desired modularity since the track and the bridge can be modeled independently. For a detailed description of this model, see Fedorova \cite{FedorovaThesis2023}.

The numerical examples presented in this section are used to motivate algorithmic choices throughout the paper. In the next section, after presenting an overview of VTSI equations, we describe encountered numerical challenges using the simplified VTSI model.

\section{Summary of VTSI governing equations and numerical challenges} \label{sec:problemStatement}

\subsection{Equations of motion}    \label{sec:eqOfMotion}

In this section, the simplified vehicle (Section \ref{sec:trainSimple}) and the simplified bridge (Section \ref{sec:bridgeSimple}) are considered. If a realistic VTSI model is considered, the equations of motion \eqref{eq:eqofmotion} are appended with equations of motion of the track and kinematic constraints between the track and the bridge. Combining the equations of motion of the train \eqref{eq:t:eqofmotion} and the bridge with the kinematic constraint equation, we obtain the coupled equations of motion

\begin{subequations} \label{eq:eqofmotion}
\begin{align}
\Mt \utddotB + \Ct(t) \utdotB + \Kt(t) \utB + {\Ltr}^{\T} \lambdaB &= \Pt(t) \label{eq:eqofmotion:t} \\ 
\Mb \ubddot + \Cb \ubdot + \Kb \ub + {\Lb (t)}^{\T} \lambdaB &= \Pb \label{eq:eqofmotion:b} \\ 
\Ltr \utB + \Lb (t) \ub &= \mathbf{0} \label{eq:eqofmotion:constr}
\end{align}
\end{subequations}
where superscripts "t", and "b" stand for train and bridge; ${\bf M^\square}$, ${\bf C^\square}$ and ${\bf K^\square}$ are the mass, damping, and stiffness matrices; ${\bf P^\square}$ is the external load vector; ${\uB^\square}$,  ${\uBdot^\square}$ and  ${\uBddot^\square}$  are vectors
of displacements, velocities, and accelerations; $\lambdaB$ is the vector of contact forces between the train wheels and the bridge. Matrices \Ltr\ and $\Lb(t)$ are kinematic constraint matrices, where $\Lb(t)$ is a time-dependent matrix derived based on the current locations of train wheels on the bridge. Matrices $\Ct(t)$ and $\Kt(t)$ are time-varying due to the corotational formulation. The vector $\lambdaB$ is the vector of Lagrange multipliers, which are widely used to enforce constraints in contact problems (\cite{Kothari2022,Duong2017,Meier2016,Kadapa2016,Matzen2016,Zhang2011,Simeon2006,Baumgarte1972}). 

System \eqref{eq:eqofmotion} contains the constraints on displacement level and is an index-3 system of differential-algebraic equations (DAE) \cite{brenan1996numerical}. The index of a DAE system is the minimum number of differentiations needed to obtain ODEs for all the unknowns. To obtain ODEs for the term ${\bf \Lambda} (t) = {\Lb (t)}^{\T} \lambdaB$, equation \eqref{eq:eqofmotion:b} has to be differentiated once. A term $\ubdddot$ then appears in the differentiated equation. Therefore, equation \eqref{eq:eqofmotion:constr} has to be differentiated three times to solve for $\ubdddot$. If the constraint equation \eqref{eq:eqofmotion:constr} is differentiated once or twice, the system can be reduced to index-2 or index-1 DAE correspondingly. We next discuss numerical challenges associated with solving an index-3 DAE system and with using straight finite elements to model the bridge.

\subsection{Numerical challenges} \label{sec:numericalChallenges}

\begin{figure}[H]
\centering
\captionsetup[subfigure]{justification=centering}
\begin{tabular}{l@{\hskip 0.8cm}r}
\subfloat[Vertical bridge acceleration] 
    {\label{fig:FEM_GenAl:a} \includegraphics{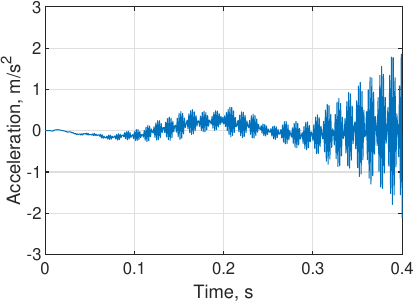}} &
\subfloat[Vertical train accelerations] 
    {\label{fig:FEM_GenAl:b} \includegraphics{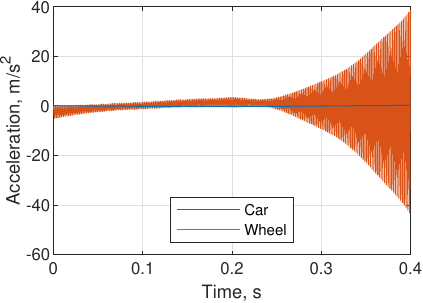}}\\
\subfloat[Vertical bridge acceleration] 
    {\label{fig:FEM_GenAl:c} \includegraphics{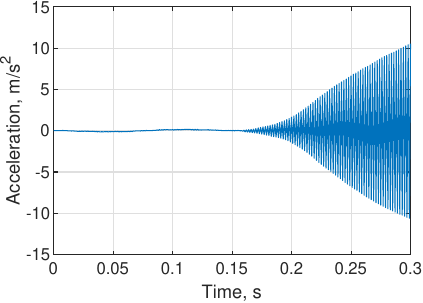}} &
\subfloat[Vertical train accelerations] 
    {\label{fig:FEM_GenAl:d} \includegraphics{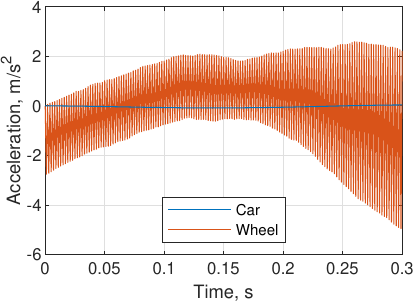}}\\
\subfloat[Vertical bridge acceleration] 
    {\label{fig:FEM_GenAl:e} \includegraphics{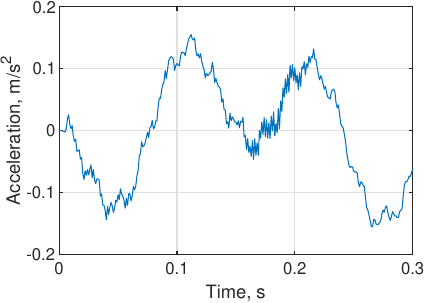}} &
\subfloat[Vertical train accelerations] 
    {\label{fig:FEM_GenAl:f} \includegraphics{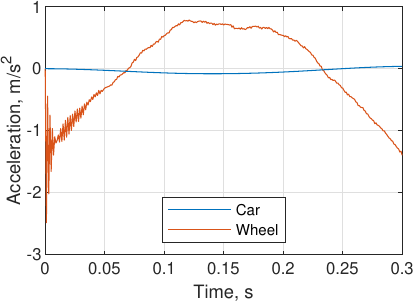}}
\end{tabular}
\caption{
    \textbf{\textit{Model:}} 
    Simplified vehicle passing with speed $100\ \si{m/s}$ over a single-span straight bridge modeled with finite elements. 
    \textbf{\textit{Algorithmic feature tested:}} 
    (a) and (b): Newmark scheme;
    (c) and (d): The Generalized-$\alpha$ scheme without numerical damping (\rhoinf = 1)\protect\footnotemark; 
    (e) and (f): The Generalized-$\alpha$ scheme with slight numerical damping (\rhoinf = 0.9).
    \textbf{\textit{Observations and Conclusions:}} 
    The presence of kinematic constraints causes unstable behavior when using numerical schemes without numerical damping. The response can be stabilized by introducing slight numerical damping.
}
\label{fig:FEM_GenAl}
\end{figure}

\begin{figure}[H]
\centering
\setcounter{subfigure}{0}
\captionsetup[subfigure]{justification=centering}
\begin{tabular}{l@{\hskip 0.8cm}r}
\subfloat[Vertical bridge acceleration (Eurocode EN 1990 (2002)\cite{EN1990} limits the value of vertical bridge deck acceleration to $3.5\ \si{m/s^2}$ for ballasted tracks)] 
    {\label{fig:FEM_GenAl09:a} \includegraphics{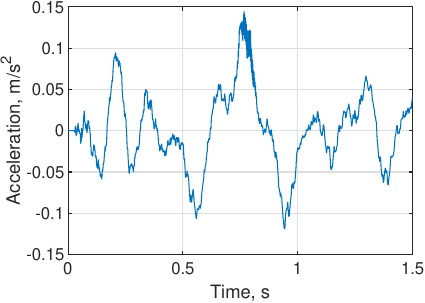}} &
\subfloat[Vertical train accelerations] 
    {\label{fig:FEM_GenAl09:b} \includegraphics{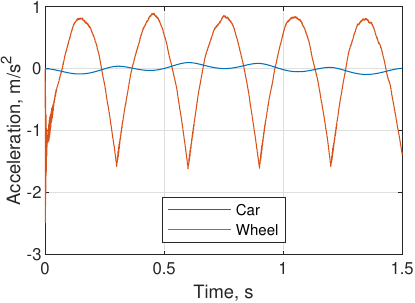}}\\
\subfloat[Transverse contact force] 
    {\label{fig:FEM_GenAl09:c} \includegraphics{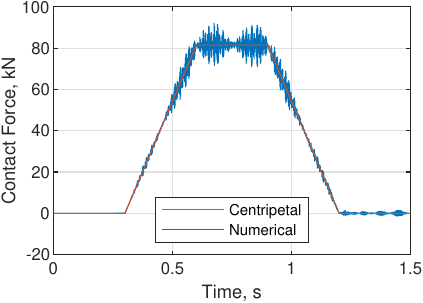}} &
\subfloat[Transverse wheel acceleration] 
    {\label{fig:FEM_GenAl09:d} \includegraphics{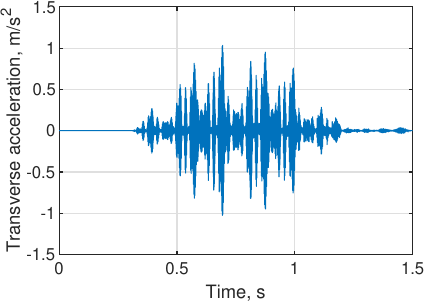}}
\end{tabular}
\caption{
    \textbf{\textit{Model:}} 
    Simplified vehicle passing with speed $100\ \si{m/s}$ over a 5-span bridge (Figure \ref{fig:modelPlanSimple}) modeled with finite elements. \textbf{\textit{Algorithmic feature tested:}} 
    The Generalized-$\alpha$ scheme (Section \ref{sec:genAlpha}) with slight numerical damping ($\rhoinf = 0.9$). \textbf{\textit{Observations and Conclusions:}} 
    Oscillations in the contact force (c) and wheel accelerations (d) in the transverse direction are amplified due to the usage of straight FEM elements in the bridge model.
}
\label{fig:FEM_GenAl09_5spans}
\end{figure}

The solution of DAE systems usually requires careful consideration (\cite{Matthies2006,Arnold2006,Borri2006,Betsch2005,Fisette1996}). When an index-3 DAE such as \eqref{eq:eqofmotion} is discretized using a time integration scheme without numerical damping, spurious oscillations are observed in wheel and bridge accelerations while using FEM to model the bridge (Figures \ref{fig:FEM_GenAl:a}-\ref{fig:FEM_GenAl:d}). This behavior, as well as oscillations in the Lagrange multipliers, is due to the presence of the kinematic constraints. As explained by Geradin et al.\cite{geradin2001flexible} and Cardona et al.\cite{Cardona1989}, a weak instability in the second time-derivative is observed when a standard Newmark scheme is used to discretize a DAE system. 
There are two strategies to mitigate these oscillations: 1) using numerical damping and 2) enforcing constraints at all three levels (displacement, velocity, and acceleration), either explicitly or through the correction steps (Section \ref{sec:constrCorrestion}). The second approach is not feasible in the case of the finite elements bridge model since standard finite elements lack inter-element curvature continuity, and, as a result, the acceleration constraint \eqref{eq:accConstr} will be discontinuous. If the first approach is employed, the introduction of numerical damping mitigates spurious oscillations in the obtained response (Figures \ref{fig:FEM_GenAl:c}-\ref{fig:FEM_GenAl:d}), similar to the effect of the Bathe discretization scheme \cite{bathe2005} in Reference \cite{Fedorova2017}. While this approach works well for mitigating the oscillations in the vertical direction, when a curved bridge is introduced, the oscillations in the transverse direction are amplified if straight FEM elements are used to approximate the curved shape of the bridge model (Figure \ref{fig:FEM_GenAl09_5spans}). Therefore, NURBS-based discretization of the bridge is introduced in Section \ref{sec:nurbs}. Such discretization provides sufficient inter-element continuity and also opens the path to enforcing acceleration constraints explicitly (Section \ref{sec:accConstr}) and correcting wheel accelerations through a simple constraint projection method (Section \ref{sec:constrCorrestion}).
%
\stepcounter{footnote}

\footnotetext{\rhoinf\ is the user-specified value of the spectral radius in the high-frequency limit. The value $\rhoinf = 1$ corresponds to no dissipation, and the value $\rhoinf = 0$ results in removing high-frequency response after the first time step (See Chung \cite{Chung1994} and Section \ref{sec:genAlpha}).}


\newpage
\section{NURBS-based beam formulation} \label{sec:nurbs}

In our past work on a two-dimensional VTSI algorithm \cite{Fedorova2017}, the Bathe method was implemented to mitigate spurious oscillations caused by the lack of $C^2$ (curvature) continuity across beam elements in the bridge model. As explained in Reference \cite{FedorovaThesis2017}, the VTSI algorithm requires the second derivative of the constraint matrix \Lb\ to be twice continuously differentiable spatially, for acceleration to be continuous in time. However, when the bridge is modeled using standard finite elements, only $C^1$ continuity in the nodes is provided.

In the present work, to provide sufficient inter-element continuity, an isogeometric approach is employed for modeling the deck in the simplified model and rails in the realistic model. This approach uses Non-Uniform Rational B-Splines (NURBS) as a basis to describe the geometry and the solution space \cite{Hughes2005}. Implementation of NURBS provides the desired $C^2$ (or higher) continuity throughout the curve, as will be shown in Section \ref{sec:algorithm}. In this study, the bridge is modeled using the three-dimensional isogeometric beam formulation by Zhang et al.~\cite{Zhang2016}. A thorough treatment of isogeometric analysis in general can be found in Hughes et al.\cite{Hughes2005} and Cottrell et al.\cite{Cottrell2009}. Algorithms for practical NURBS implementation are provided in Piegl et al. \cite{Piegl1997}. For the use of NURBS in contact mechanics, see, for example, Agrawal et al. \cite{Agrawal2020}. Here, we provide a brief explanation of the modeling technique. 

\subsection{NURBS-based representation of geometry and solution fields}
 
A NURBS curve can then be defined as (see Cottrell et al. \cite{Cottrell2009} and Piegl et al. \cite{Piegl1997})
\begin{equation} \label{eq:nurbs:curve}
    \xtilde(\xi) = \sum\limits_{i = 1}^{n} R_{i,p}(\xi) \x_i
\end{equation}
where $\x_i$ is the $i$-th control point, $n$ is the number of basis functions, and $R_{i,p}(\xi)$ are the rational basis functions defined as
\begin{equation} \label{eq:nurbs:shapeFunctions}
    R_{i,p}(\xi) = \frac{N_{i,p}(\xi) w_i}{\sum \limits_{j = 1}^{n} N_{j,p}(\xi) w_j} = \frac{N_{i,p}(\xi) w_i}{W(\xi)}
\end{equation}
where $w_i$ are the control points' weights and $N_{i,p}(\xi)$, $i = 1,2,...,n$, are the B-spline basis functions of degree $p$ defined on a knot vector 
\begin{equation}    \label{eq:nurbs:knotVector}
    \Xi = \{ \underbrace{0,...,0}_{p+1} , 1, 2, ..., n^e-1,
    \underbrace{n^e,...,n^e}_{p+1} \} 
\end{equation}
where the number of elements $n^e = n-p$ is equal to the number of nonzero knot spans.

Equation \ref{eq:nurbs:curve} can be written in a matrix form as
\begin{equation} \label{eq:b:xhat}
    \xtilde(\xi) = \mathbf{R}^\T(\xi) \x  
\end{equation} 

The geometry of the curve can either be represented exactly or, given a set of data points, it can be approximated, for example, using least square approximation. In either case, the NURBS basis is used to describe both geometry and solution fields. Displacement and rotation fields can then be approximated similarly to \eqref{eq:b:xhat} as
\begin{equation}
\begin{gathered}
    \upsilontilde(\xi) = \sum\limits_{i = 1}^{n} R_{i,p}(\xi) \upsilonb_i \\
    \thetatilde(\xi) = \sum\limits_{i = 1}^{n} R_{i,p}(\xi) \thetab_i
\end{gathered}
\end{equation}
where \upsilonb\ and \thetab\ are control displacements and rotations. 

Knots of the knot vector, $\Xi$, mesh the curve $\xtilde$ into elements. Each element $\Omegae = [ \xi_i,\xi_{i+1} ]$ is supported by $p+1$ basis functions $R_{i-p,p},...,R_{i,p}$. The geometry, displacement, and rotation fields at the point $\xi \in \Omegae$ can then be expressed as \cite{Zhang2016}
\begin{equation}
\begin{gathered}
    \xtilde(\xi) = \sum\limits_{k = i-p}^{i} R_{k,p}(\xi) \x_k = 
    {\Rmem}^\T \xe\\
    \upsilontilde(\xi) = \sum\limits_{k = i-p}^{i} R_{k,p}(\xi) \upsilonb_k = {\Rmem}^\T \upsilonbe\\
    \thetatilde(\xi) = \sum\limits_{k = i-p}^{i} R_{k,p}(\xi) \thetab_k = 
    {\Rmem}^\T \thetabe\\
\end{gathered}
\end{equation}
%

\subsection{Kinematics of the Timoshenko beam}

Timoshenko beam theory is employed. Therefore, the generalized strains can be expressed as \cite{Zhang2016, Tabarrok1988}
\begin{equation} \label{eq:strain}
\begin{gathered}
    \boldsymbol{\Theta}(s) = \left( {\begin{array}{*{20}{c}}
    \boldsymbol{\varepsilon}(s) \\
    \boldsymbol{\beta}(s) 
    \end{array}} \right) = 
    \left( {\begin{array}{c}
    \frac{d\boldsymbol{\upsilon}}{ds} - \boldsymbol{\theta} \times\textbf{t}\\
    \frac{d\boldsymbol{\theta}}{ds}    
    \end{array}} \right) = 
    \left( {\begin{array}{l}
    \frac{d\upsilon_t}{ds} - \kappa \upsilon_n\\[0.1cm]
    \frac{d\upsilon_n}{ds} + \kappa \upsilon_t - \tau \upsilon_b - \theta_b\\[0.1cm]
    \frac{d\upsilon_b}{ds} + \tau \upsilon_n + \theta_n\\[0.1cm]
    \frac{d\theta_t}{ds} - \kappa \theta_n\\[0.1cm]
    \frac{d\theta_n}{ds} + \kappa \theta_t - \tau \theta_b\\[0.1cm]
    \frac{d\theta_b}{ds} - \tau \theta_n
    \end{array}} \right)   
\end{gathered}
\end{equation}

\noindent where $\boldsymbol{\varepsilon}$ is the vector of axial and the two shear strains, $\boldsymbol{\beta}$ is the vector of one twisting and the two bending strains,
$\boldsymbol{\upsilon}$ and $\boldsymbol{\theta}$ are independent displacement and rotation fields in Frenet basis $\{\textbf{t},\textbf{n},\textbf{b}\}$, and $s$ is the arc-length parameter. The parameter $\xi$ is mapped to the arc length as $\xi = \hat{\xi}(s)$, as described in Section \ref{sec:corotational}. Derivatives $\frac{d(\ )}{ds}$ are found through the chain rule as $\frac{d(\ )}{ds} = \frac{d(\ )}{\xi} \frac{d\xi}{ds} = \frac{d(\ )}{\xi} \frac{1}{J}$, where the Jacobian is computed as $J = \frac{ds}{d\xi} = \sqrt{(\frac{dx}{d\xi})^2 + (\frac{dy}{d\xi})^2 + (\frac{dz}{d\xi})^2}$.

\subsection{Energy of the Timoshenko beam}

The strain energy of the Timoshenko beam can be expressed as
\begin{equation} \label{eq:strainEnergy}
    U = \frac{1}{2} \int\limits_{0}^{L}
    \left(
    EA \varepsilon_t^2 +
    GA_n \varepsilon_n^2 +
    GA_b \varepsilon_b^2 +
    GI_t \beta_t^2 +
    EI_n \beta_n^2 +
    EI_b \beta_b^2
    \right) ds
\end{equation}

Having the generalized strain expressions \eqref{eq:strain} substituted into equation \eqref{eq:strainEnergy}, derivatives of \eqref{eq:strainEnergy} with respect to displacement and rotation fields are found.

The kinetic energy of the beam is defined as
\begin{equation} \label{eq:b:kinetic}
    T = \frac{1}{2} \int\limits_{0}^{L} 
    \left( \begin{array}{l}
         \upsilondot  \\
         \thetadot
    \end{array} \right)^T
    \rho 
    \left[ \begin{array}{ll}
         A \textbf{I} & {}\\
         {} & \textbf{I}^{\text{b}}
    \end{array} \right]
    \left( \begin{array}{l}
         \upsilondot  \\
         \thetadot
    \end{array} \right) ds
\end{equation}
where $A$ is the cross-sectional area, $\textbf{I}$ is an identity matrix, and $\textbf{I}^{\text{b}}$ is a diagonal matrix with components $\{I_t,I_n,I_b\}$, which are area moments of inertia about $t$, $n$, and $b$ axes, correspondingly.

Work done by distributed loads, such as self-weight, can be found as
\begin{equation} \label{eq:b:work}
    V = \frac{1}{2} \int\limits_{0}^{L} 
    \rho g \boldsymbol{\upsilon}_b ds
\end{equation}
The element load vector can then be found. Further details on the derivation of member matrices can be found in Zhang et al.~\cite{Zhang2016} and Fedorova~\cite{FedorovaThesis2023}.

\section{Time integration} \label{sec:algorithm}

In this Section, three ways of solving the VTSI equations are described. The three methods have their advantages and disadvantages. Some are easy to implement in stand-alone codes. Others are better suited to be incorporated into an existing structural analysis software.

\subsection{Generalized-$\alpha$ method with numerical dissipation} \label{sec:genAlpha}

As noted earlier, the system of VTSI equations \eqref{eq:eqofmotion} is an index-3 DAE system. As the available research suggests, some amount of numerical damping is necessary to avoid spurious high-frequency numerical oscillations in the Lagrange multipliers in such systems \cite{Arnold2007}. Therefore, the Generalized-$\alpha$ method \cite{Chung1994} is employed to discretize equations of motion in time. This scheme provides accuracy at low-frequency and desired numerical damping at high-frequency \cite{Arnold2007}. Other methods, such as Backward Difference Formulas and the Bathe method \cite{bathe2005} could also be used to this end; however, the Generalized-$\alpha$ method is consistent with implementation in several existing bridge analysis software, making it more convenient to integrate the proposed VTSI algorithm, and is often used in discretization of constrained systems (\cite{Bruls2008,Bruls2014}).

The VTSI equations of motion \eqref{eq:eqofmotion:t}-\eqref{eq:eqofmotion:b} are discretized in time using averaging parameters \alm\ and \alf\ of the scheme as follows 

\begin{subequations} \label{eq:discr}
\begin{align}
\Mt \utddotB_{\nalm} + \Ct(\tnalf) \utdotB_{\nalf}  + \Kt(\tnalf) \utB_{\nalf}  + {\Ltr}^\T \lambdaB_{\nalf}  &= \Pt (\tnalf) \label{eq:discr:a} \\ 
\Mb \ubddot_{\nalm}  + \Cb \ubdot_{\nalf}  + \Kb \ub_{\nalf}  + \left( \Lb(\tnalf) \right)^\T \lambdaB_{\nalf}  &= \Pb \label{eq:discr:b}
\end{align}
\end{subequations}

Even though Generalized-$\alpha$ is not an energy-preserving scheme, in order to discretize the constraints, we assume that the work done by the constraint forces over the time increment should vanish \cite{geradin2001flexible,bauchau2010flexible} (as for an energy-preserving scheme). This can be written using Equation \eqref{eq:euler2} as

\begin{equation}    \label{eq:constrWork}
    W \approx (u_{i,\ n+1} - u_{i,\ n}) \left( \frac{\partial C_j}{\partial u_i} \right)_{\nalf} \lambda_{j,\ \nalf} 
    \approx \left( C_{j,\ \nplusone} - C_{j,\ n} \right) \lambda_{j,\ \nalf}
\end{equation}
%
%



It follows from Equation \eqref{eq:constrWork} that enforcing constraint ${\left(C_j \right) }\nplusone$ will ensure the vanishing of the work done by constraint forces. Therefore, constraint equation \eqref{eq:eqofmotion:constr} is discretized in time as 
\begin{equation} \label{eq:discr:constr}
\Ltr \utB_{\nplusone}  + \Lb(t_{n+1}) \ub_{\nplusone} = \mathbf{0}
\end{equation}


The analysis steps are outlined in Procedure \ref{procedure_genAlpha}. The details of the time discretization of Equations \eqref{eq:discr}, as well as the notation of Procedure~\ref{procedure_genAlpha}, can be found in Fedorova \cite{FedorovaThesis2023}. Steps particular to 3D analysis and curved paths are highlighted.

\begin{figure}[H]
\centering
\setcounter{subfigure}{0}
\captionsetup[subfigure]{justification=centering}
\begin{tabular}{l@{\hskip 0.8cm}r}
\subfloat[Vertical bridge acceleration] 
    {\label{fig:NURBS_GenAl09:a} 
    \includegraphics{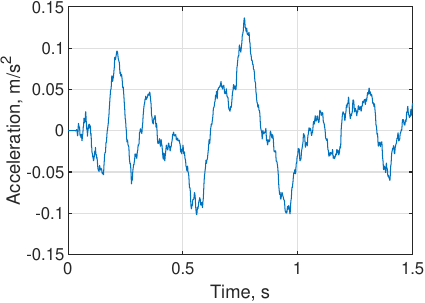}} &
\subfloat[Vertical train accelerations] 
    {\label{fig:NURBS_GenAl09:b} \includegraphics{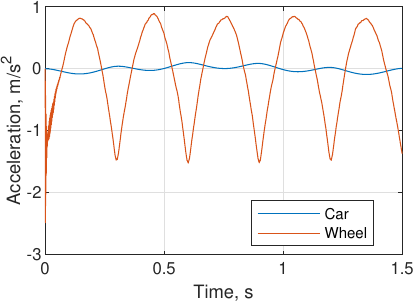}}\\
\subfloat[Transverse contact force] 
    {\label{fig:NURBS_GenAl09:c} 
    \includegraphics{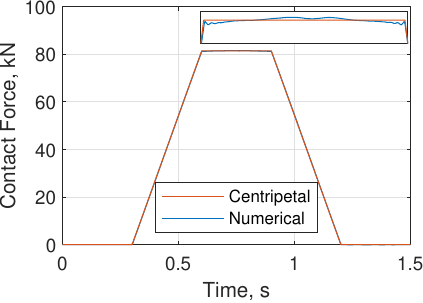}} &
\subfloat[Transverse wheel acceleration] 
    {\label{fig:NURBS_GenAl09:d} \includegraphics{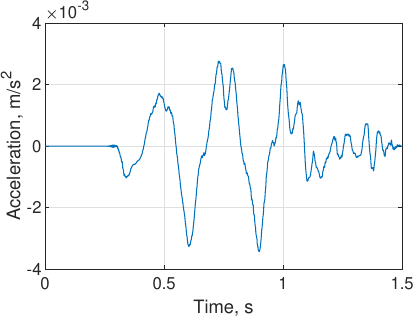}}
\end{tabular}
\caption{
    \textbf{\textit{Model:}} 
    Simplified vehicle passing with speed $100\ \si{m/s}$ over a 5-span bridge (Figure \ref{fig:modelPlanSimple}) modeled with NURBS of degree $p = 3$. 
    \textbf{\textit{Algorithmic feature tested:}} 
    Generalized-$\alpha$ scheme (Section \ref{sec:genAlpha}) with slight numerical damping ($\rhoinf = 0.9$). 
    \textbf{\textit{Observations and Conclusions:}} 
    c) implementation of NURBS helps to mitigate spurious oscillations in the transverse direction (compare with Figures \ref{fig:FEM_GenAl09:c} and \ref{fig:FEM_GenAl09:d}, where FEM elements are used. \textit{Note}: The difference between the numerical result and the theoretical is due to the fact that centripetal force is calculated for the mass element only, not accounting for the effect of the suspension element);
    b) acceleration discontinuity at time $t = 0$: the vehicle enters the bridge, and the curvature changes abruptly.
}
\label{fig:NURBS_GenAl09}
\end{figure}

Figure \ref{fig:NURBS_GenAl09} illustrates the results of VTSI analysis of the 5-span simplified bridge (Figure \ref{fig:modelPlanSimple}) modeled with NURBS and the simplified vehicle model (Figure \ref{fig:modelSimple}). As can be seen in Figure \ref{fig:NURBS_GenAl09:c}, implementation of NURBS helps to mitigate spurious oscillations in the transverse direction (compare with Figure \ref{fig:FEM_GenAl09:c}, where FEM elements are used). However, Figure \ref{fig:NURBS_GenAl09:b} clearly demonstrates the presence of acceleration discontinuity at time $t = 0$. At this time, the vehicle enters the bridge, and the curvature changes abruptly. The next section discusses ways to avoid such discontinuity.

\makeatletter
\xpatchcmd{\algorithmic}{\itemsep\z@}{\itemsep=1ex}{}{}
\makeatother

\floatname{algorithm}{Procedure} 
\begin{algorithm} 
\caption{\text Linear dynamic analysis of Vehicle-Track-Structure Interaction} \label{procedure_genAlpha}
\begin{algorithmic}[1]
\State Given: \Pt, \Pb, $\rhoB (t)$, \dt, \alf, \alm, $\beta$, $\gamma$, bridge geometry, train geometry
\State Assemble \Mt, \Ltr
\State Assemble \Mb, \Cb, \Kb
\State Assemble \Mbbar
\State Compute $\utB_0$ by linear static analysis. Set $\utdotB_0 = \bf{0}$ \Comment{Train initial conditions}
\State Compute $\ub_0$ by linear static analysis. Set $\ubdot_0 = \bf{0}$ \Comment{Bridge initial conditions}
\For{$n\leftarrow 0$, NUMINC} 
\State Compute $\tnplus = t_n + \dt$, $\tnalf = (1-\alf) t_{n+1} + \alf t_n$
\State Compute the position of the wheel at the time step \tnplus , \tnalf
\State \tikzmk{A}Compute Frenet frame properties $\RF\nalf$, $\xFddot\nalf$, $\omegaFhat\nalf$, $\omegaFhatdot\nalf$ \tikzmk{B} \boxit{gray}
\State Obtain $\Lb ( \tnplus )$, $\Lb ( \tnalf )$
\State Compute $\ab\nalf$ 
\State \tikzmk{A}Compute $\Pt\nalf$, $\Mtbar\nalf$, $\at\nalf$
    \tikzmk{B} \boxit{gray}
\State Solve $\Mtbar \left[\uttilde \quad \At \right]= \left[\at \quad {\Ltr}^\T \right]$
\State Solve $\Mbbar \left[\ubtilde \quad \Ab \right]= \left[\ab \quad {\left( \Lb (\tnalf) \right)}^{\T} \right]$
\State Set $\Abar =  \left[ \Ltr \At + \Lb (\tnplus) \Ab \right]$
\State Set $\rhoBbar =  \Ltr \uttilde + \Lb (t_{n+1}) \ubtilde$
\State Solve $\Abar \lambdaB\nalf = \rhoBbar$
\State Set $\utB_{n+1} = \uttilde - \At \lambdaB{\nalf}$; \quad
$\utddotB_{n+1} = \frac{1}{\dt^2 \beta} \left( \utB_{n+1} - \utB_n - \dt \utdotB_n \right) - \left( \frac{1}{2 \beta} - 1 \right) \utddotB$; \hfill\break \hspace*{1.2cm}
$ \utdotB_{n+1} = \utdotB_n + \dt \left( (1-\gamma) \utddotB_n + \gamma \utddotB_{n+1} \right) $
\State Set $\ub_{n+1} = \ubtilde - \Ab \lambdaB{\nalf}$; \quad
$\ubddot_{n+1} = \frac{1}{\dt^2 \beta} \left( \ub_{n+1} - \ub_n - \dt \ubdot_n \right) - \left( \frac{1}{2 \beta} - 1 \right) \ubddot$; \hfill\break \hspace*{1.2cm}
$\ubdot_{n+1} = \ubdot_n + \dt \left( (1-\gamma) \ubddot_n + \gamma \ubddot_{n+1} \right) $

\EndFor
\end{algorithmic}
\end{algorithm}

\subsection{Explicit enforcement of acceleration constraints} \label{sec:accConstr}

One way to avoid the accumulation of error in obtained accelerations (see Figures \ref{fig:FEM_GenAl:a}-\ref{fig:FEM_GenAl:b}) is to enforce the acceleration constraint \eqref{eq:accConstr} explicitly instead of the displacement constraint \eqref{eq:discr:constr}. Such an approach is traditionally used in multibody dynamics (for example, see Shabana \cite{shabana2001computational}). In this case, the VTSI equations of motion are the same as \eqref{eq:eqofmotion:t}-\eqref{eq:eqofmotion:b}, and the acceleration constraint can be discretized similarly to \eqref{eq:discr:constr} as

\begin{equation} \label{eq:discr:accConstr}
\Ltr \utddotB_{\nplusone} + \ddLb(\tnplus) \ub_{\nplusone} + 2\dLb(\tnplus) \ubdot_{\nplusone} + \Lb(\tnplus) \ubddot_{\nplusone} = 0 
\end{equation}
The system of Equations \eqref{eq:discr} and \eqref{eq:discr:accConstr} can be rewritten similarly to how the original system of Equations \eqref{eq:discr} and \eqref{eq:discr:constr} was rewritten for the Procedure \ref{procedure_genAlpha}.


\begin{figure}[H]
\centering
\setcounter{subfigure}{0}
\captionsetup[subfigure]{justification=centering}
\begin{tabular}{l@{\hskip 0.8cm}r}
\subfloat[Vertical bridge acceleration] 
    {\label{fig:NURBS_Newmark_accConstr:p3:a} 
    \includegraphics{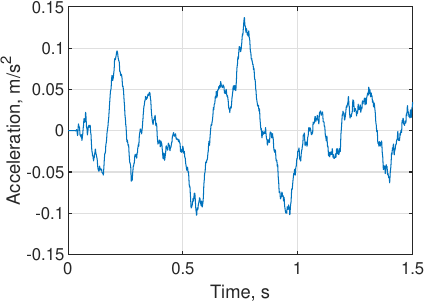}} &
\subfloat[Vertical car accelerations] 
    {\label{fig:NURBS_Newmark_accConstr:p3:b} \includegraphics{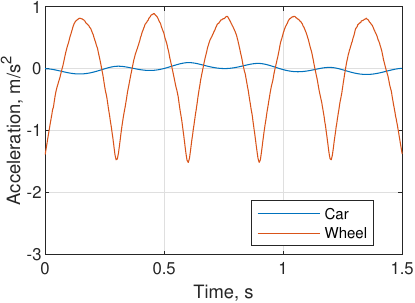}}\\
\subfloat[Transverse contact force] 
    {\label{fig:NURBS_Newmark_accConstr:p3:c} 
    \includegraphics{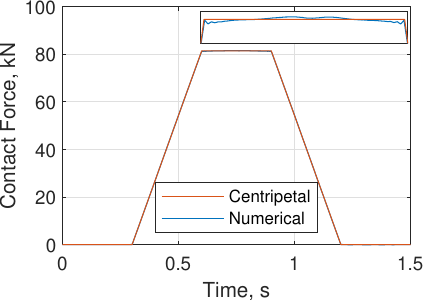}} &
\subfloat[Transverse wheel accelerations] 
    {\label{fig:NURBS_Newmark_accConstr:p3:d} \includegraphics{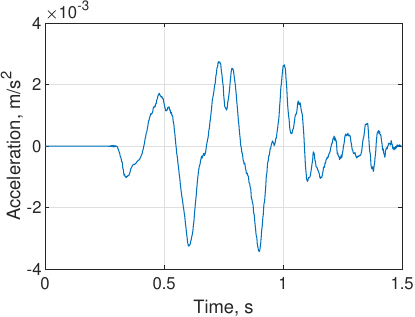}}
\end{tabular}
\caption{
    \textbf{\textit{Model:}} 
    Simplified vehicle passing with speed $100\ \si{m/s}$ over a 5-span bridge (Figure \ref{fig:modelPlanSimple}) modeled with NURBS of \textbf{degree \textit{p} = 3}. 
    \textbf{\textit{Algorithmic feature tested:}} 
    Acceleration constraints are enforced explicitly (Equation \eqref{eq:discr:accConstr}) at each time step. Wheel velocity and acceleration are corrected only at $t=0$ using constraints \eqref{eq:correction}. Newmark scheme ($\alm = 0,\ \alf = 0$) is used.
    \textbf{\textit{Observations and Conclusions:}} 
    Enforcement of acceleration constraint instead of displacement constraint allows eliminating the spurious oscillations without employing the Generalized-$\alpha$ scheme and by correcting only initial values of wheel velocity and acceleration.
}
\label{fig:NURBS_Newmark_accConstr_p3}
\end{figure}

\begin{figure}[H]
\centering
\setcounter{subfigure}{0}
\captionsetup[subfigure]{justification=centering}
\begin{tabular}{l@{\hskip 0.8cm}r}
\subfloat[Vertical bridge acceleration] 
    {\label{fig:NURBS_Newmark_accConstr:p5:a} 
    \includegraphics{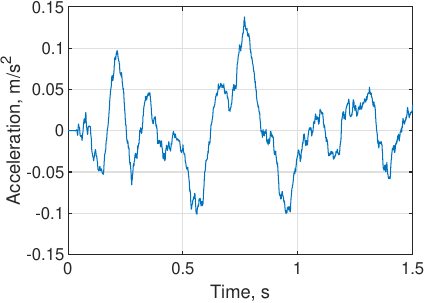}} &
\subfloat[Vertical car accelerations] 
    {\label{fig:NURBS_Newmark_accConstr:p5:b} \includegraphics{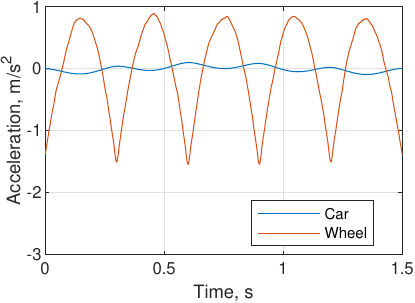}}\\
\subfloat[Transverse contact force] 
    {\label{fig:NURBS_Newmark_accConstr:p5:c} 
    \includegraphics{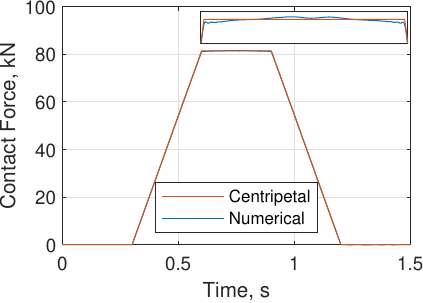}} &
\subfloat[Transverse wheel accelerations] 
    {\label{fig:NURBS_Newmark_accConstr:p5:d} \includegraphics{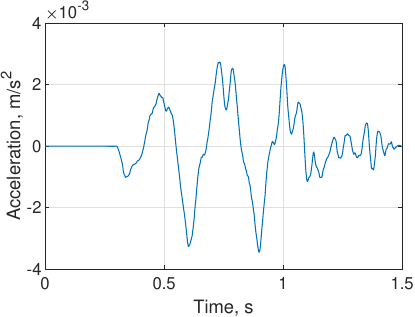}}
\end{tabular}
\caption{
    \textbf{\textit{Model:}} 
    Simplified vehicle passing with speed $100\ \si{m/s}$ over a 5-span bridge (Figure \ref{fig:modelPlanSimple}) modeled with NURBS \textbf{degree \textit{p} = 5}. 
    \textbf{\textit{Algorithmic feature tested:}} 
    Acceleration constraints are enforced explicitly (Equation \eqref{eq:discr:accConstr}) at each time step. Wheel velocity and acceleration are corrected only at $t=0$ using constraints \eqref{eq:correction}. Newmark scheme ($\alm = 0,\ \alf = 0$) is used.
    \textbf{\textit{Observations and Conclusions:}} 
    Enforcement of acceleration constraint instead of displacement constraint allows eliminating the spurious oscillations without employing the Generalized-$\alpha$ scheme and by correcting only initial values of wheel velocity and acceleration. A higher degree of NURBS helps to further mitigate the oscillations (compare with Figure \ref{fig:NURBS_Newmark_accConstr_p3}).
}
\label{fig:NURBS_Newmark_accConstr_p5}
\end{figure}

\begin{figure}[H]
\centering
\setcounter{subfigure}{0}
\captionsetup[subfigure]{justification=centering}
\begin{tabular}{l@{\hskip 0.8cm}r}
\subfloat[Displacement drift without correction] 
    {\label{fig:NURBS_Newmark_accConstr:displDrift:notCorrected} 
    \includegraphics{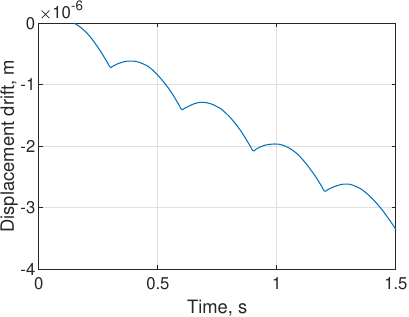}} &
\subfloat[Displacement drift with a correction of wheel displacement at $t = 0.75\ s$] 
    {\label{fig:NURBS_Newmark_accConstr:displDrift:corrected} \includegraphics{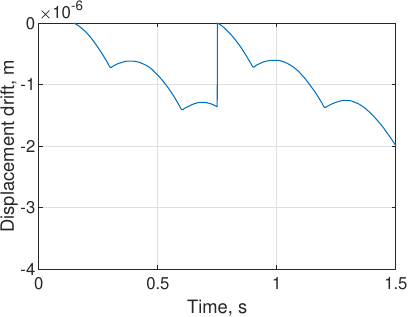}}
\end{tabular}
\caption{
    \textbf{\textit{Model:}} 
    Simplified vehicle passing with speed $100\ \si{m/s}$ over a 5-span bridge (Figure \ref{fig:modelPlanSimple}) modeled with NURBS degree $p = 5$. 
    \textbf{\textit{Observations:}} 
    Drift in the displacement constraint occurs when only acceleration constraint \eqref{eq:discr:accConstr} is enforced explicitly at each time step.
}
\label{fig:NURBS_Newmark_accConstr:displDrift}
\end{figure}

As can be seen from Figures \ref{fig:NURBS_Newmark_accConstr_p3}-\ref{fig:NURBS_Newmark_accConstr_p5}, enforcement of acceleration constraint instead of displacement constraint allows to eliminate the spurious oscillations without employing the Generalized-$\alpha$ scheme. Wheel velocity and acceleration must be corrected only at $t=0$ using constraints \eqref{eq:correction}. Enforcing only acceleration constraints throughout the analysis leads to drift in displacement constraints (see Figure \ref{fig:NURBS_Newmark_accConstr:displDrift:notCorrected}). Wheel displacements can be corrected as

\begin{equation} \label{eq:correction:displ}
    \uw_{\nplusone}  = \Lb(t_{n+1}) \ub_{\nplusone} 
\end{equation}
Figure \ref{fig:NURBS_Newmark_accConstr:displDrift:corrected} shows the result of the correction at a single time step $t = 0.75\ s$. In general, the correction of displacement constraints can be performed periodically if needed.

\subsection{Constraint projection}  \label{sec:constrCorrestion}
NURBS shape functions provide desired curvature continuity inside the bridge span. However, at time $t=0$, when the train has not entered the bridge yet, the continuity of velocity and acceleration constraints is not guaranteed.

Consider differentiating equation \eqref{eq:eqofmotion:constr} twice to obtain the constraints at velocity and acceleration level

\begin{subequations} \label{eq:constrDiffer}
    \begin{align}
        \Ltr \utdotB + \dLb\ub + \Lb\ubdot &= 0 \label{eq:velConstr}\\      
        \Ltr \utddotB + \ddLb\ub + 2\dLb\ubdot + \Lb\ubddot &= 0  \label{eq:accConstr}
    \end{align}
\end{subequations}
Recognizing that terms $\Ltr \utdotB$ and $\Ltr \utddotB$ correspond to wheel velocity and acceleration, we can write

\begin{subequations} \label{eq:correction}
    \begin{align}
        \duw &= \dLb\ub + \Lb\ubdot \label{eq:correction:a}\\      
        \dduw &= \ddLb\ub + 2\dLb\ubdot + \Lb\ubddot  \label{eq:correction:b}
    \end{align}
\end{subequations}

At time $t=0$ wheel velocity and acceleration are equal to zero. At the same time, the right-hand sides of equations \eqref{eq:correction} will not be equal to zero at this time step, resulting in the discontinuities in wheel velocity and acceleration. To overcome this, Equations \eqref{eq:correction} can be used to correct wheel velocity and acceleration at time step $t=0$.


\subsubsection{Simple example: Bridge modeled as a rigid cosine curve} 

\begin{figure}[H]
\centering
\setcounter{subfigure}{0}
\captionsetup[subfigure]{justification=centering}
\begin{tabular}{l@{\hskip 0.8cm}r}
\subfloat[Vertical train accelerations without correction] 
    {\label{fig:cos:a} \includegraphics{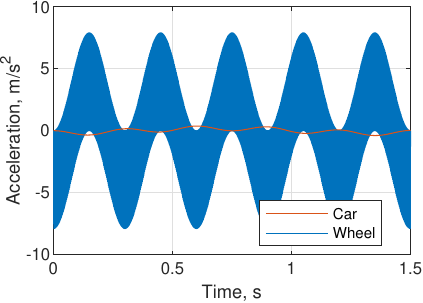}}&
\subfloat[Vertical train accelerations with a correction \eqref{eq:correction:b} of wheel acceleration at $t = 0$] 
    {\label{fig:cos:b} 
    \includegraphics{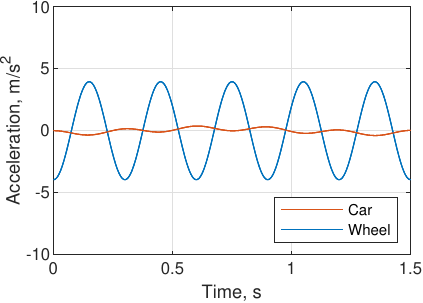}}
\end{tabular}
\caption{
    \textbf{\textit{Model:}} 
    Simulation of a vehicle running with speed $100\ \si{m/s}$ over a rigid cosine curve.
    \textbf{\textit{Algorithmic feature tested:}} 
    b) acceleration correction is performed only at $t = 0$, while displacement constraints are enforced explicitly (Equation \eqref{eq:discr:constr}) at each time step.
    \textbf{\textit{Observations and Conclusions:}} 
    By correcting the acceleration at time step $t=0$ using Equations \eqref{eq:correction:b}, we avoid discontinuities in the constraints at the acceleration level and spurious oscillations in the accelerations.
}
\label{fig:cos}
\end{figure}

Consider a simple example where the bridge in replaced with a rigid cosine curve. By correcting the acceleration at time step $t=0$ using Equations \eqref{eq:correction}, we avoid discontinuities in the constraints at the acceleration level and spurious oscillations in the accelerations (Figure \ref{fig:cos}).


\subsubsection{Example: Bridge modeled using NURBS} 

\begin{figure}[H]
\centering
\setcounter{subfigure}{0}
\captionsetup[subfigure]{justification=centering}
\begin{tabular}{l@{\hskip 0.8cm}r}
\subfloat[Vertical bridge acceleration] 
    {\label{fig:NURBS_GenAl_corr:a} 
    \includegraphics{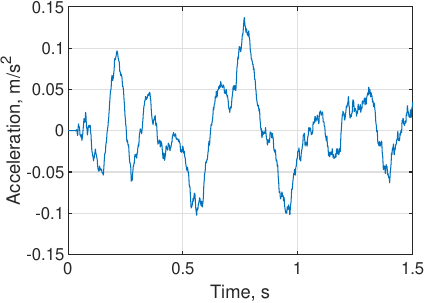}} &
\subfloat[Vertical train accelerations] 
    {\label{fig:NURBS_GenAl_corr:b} \includegraphics{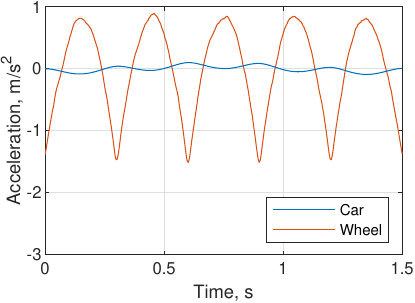}}\\
\subfloat[Transverse contact force] 
    {\label{fig:NURBS_GenAl_corr:c} 
    \includegraphics{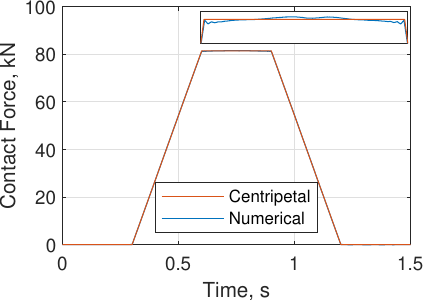}} &
\subfloat[Transverse wheel acceleration] 
    {\label{fig:NURBS_GenAl_corr:d} \includegraphics{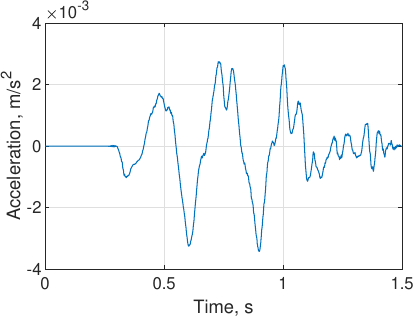}}
\end{tabular}
\caption{
    \textbf{\textit{Model:}} 
    Simplified vehicle passing with speed $100\ \si{m/s}$ over a 5-span bridge (Figure \ref{fig:modelPlanSimple}) modeled with NURBS of degree $p = 3$. 
    \textbf{\textit{Algorithmic feature tested:}} 
    Wheel velocity and acceleration corrections \eqref{eq:correction} are performed at each time step, while displacement constraints are enforced explicitly (Equation \eqref{eq:discr:constr}). Newmark scheme ($\alm = 0,\ \alf = 0$) is used.
    \textbf{\textit{Observations and Conclusions:}} 
    The correction of velocity and acceleration constraints at each time step helps to avoid the accumulation of errors in accelerations and makes the use of the Newmark scheme possible.
}
\label{fig:NURBS_GenAl_corr}
\end{figure}

In a realistic case, when displacement constraints are enforced explicitly using Equation \eqref{eq:discr:constr}, the correction only at time $t = 0$ is not sufficient since the error in accelerations accumulates and has to be corrected at each time step using Equations \eqref{eq:correction}. Figure \ref{fig:NURBS_GenAl_corr} illustrates the case when velocity and acceleration corrections are performed at each time step. 

\begin{figure}[H]
    \centering
    {\includegraphics[scale = 1.0]{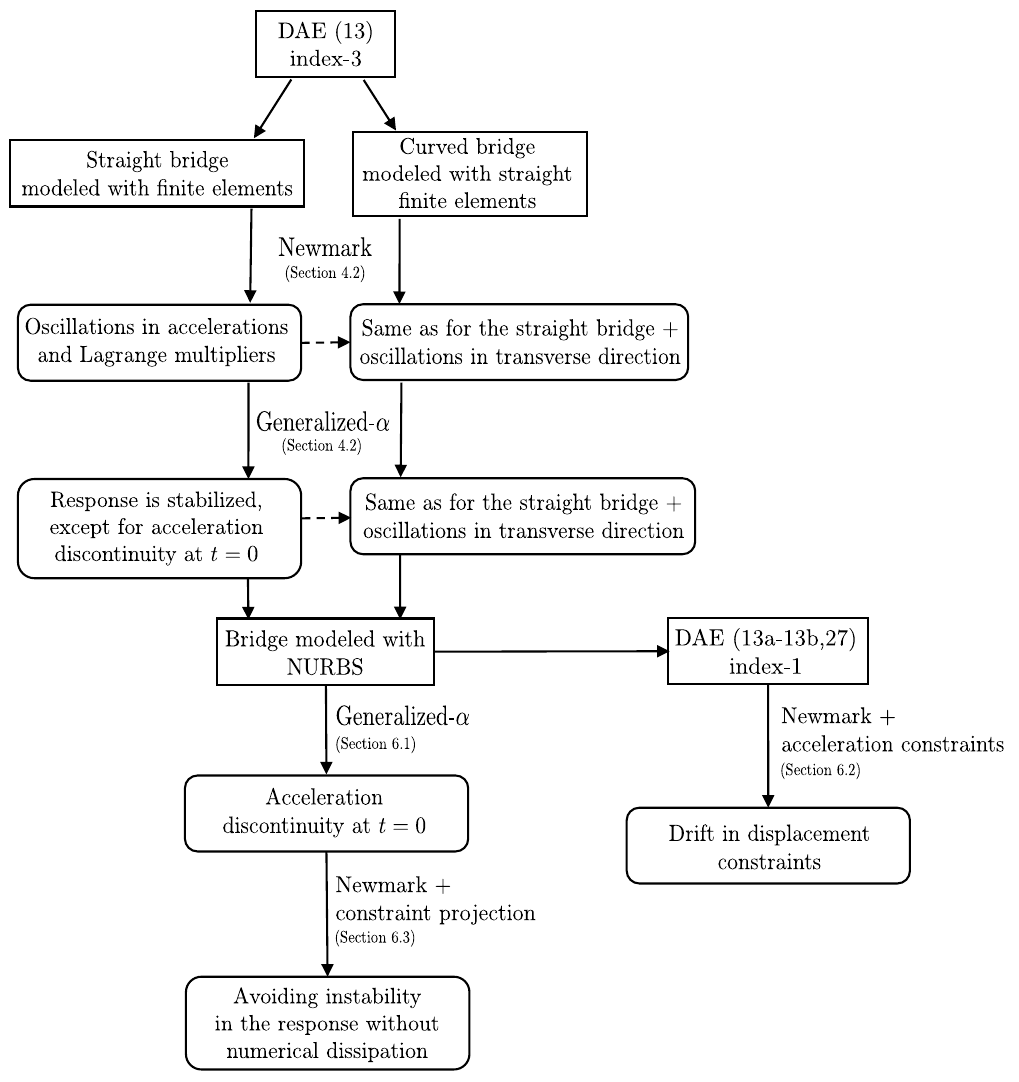}}
    \caption{Summary of algorithmic choices}
    \label{fig:summary}
\end{figure}

\section{Numerical results} \label{sec:moreNumResults}

\subsection{Simplified model}

In this section, further numerical results for the simplified VTSI model (Sections \ref{sec:trainSimple}-\ref{sec:bridgeSimple}) are presented. These results are obtained using the Generalized-$\alpha$ scheme with numerical dissipation (Section \ref{sec:genAlpha}), since among the presented choices, it is the most suitable way to integrate the VTSI algorithm into an existing structural analysis software. Vertical bridge and train accelerations, as well as transverse wheel acceleration and transverse contact force, were already presented in Figure \ref{fig:NURBS_GenAl09}. 

\begin{figure}[H]
\centering
\setcounter{subfigure}{0}
\captionsetup[subfigure]{justification=centering}
\begin{tabular}{l@{\hskip 0.8cm}r}
\subfloat[Vertical bridge displacements in the middle of the structure] 
    {\label{fig:NURBS_GenAl09_simple:a} 
    \includegraphics{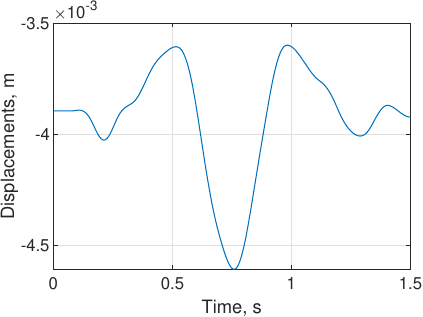}} &
\subfloat[Transverse bridge displacements in the middle of the structure] 
    {\label{fig:NURBS_GenAl09_simple:b} \includegraphics{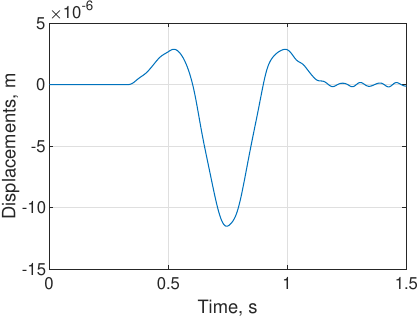}}\\
\subfloat[Vertical vehicle displacements] 
    {\label{fig:NURBS_GenAl09_simple:c} 
    \includegraphics{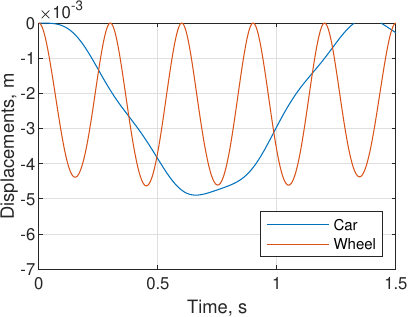}} &
\subfloat[Transverse vehicle displacements and rotations] 
    {\label{fig:NURBS_GenAl09_simple:d} \includegraphics{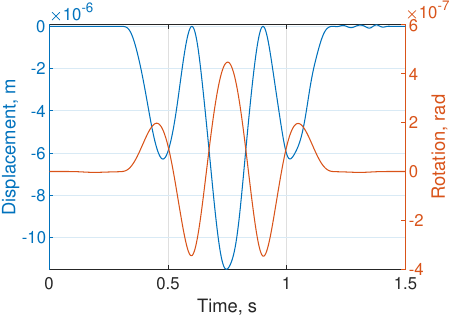}}\\
\subfloat[Vertical contact force] 
    {\label{fig:NURBS_GenAl09_simple:e} 
    \includegraphics{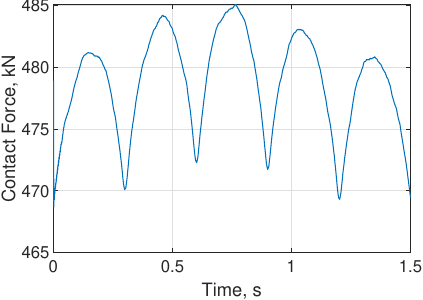}} &
\subfloat[Contact force in the direction $\theta_n$] 
    {\label{fig:NURBS_GenAl09_simple:f} \includegraphics{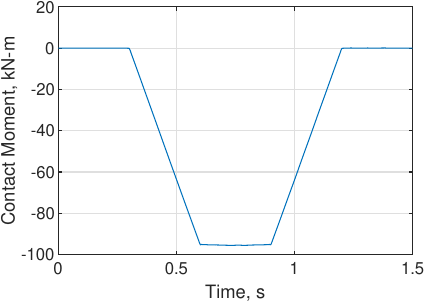}}
\end{tabular}
\caption{
    \textbf{\textit{Model:}} 
    Simplified vehicle passing with speed $100\ \si{m/s}$ over a 5-span bridge (Figure \ref{fig:modelPlanSimple}) modeled with NURBS of degree $p = 3$. 
    \textbf{\textit{Algorithmic feature tested:}} 
    The Generalized-$\alpha$ scheme (Section \ref{sec:genAlpha}) with slight numerical damping ($\rhoinf = 0.9$). 
}
\label{fig:NURBS_GenAl09_simple}
\end{figure}

As can be seen from Figures \ref{fig:NURBS_GenAl09_simple} and \ref{fig:NURBS_GenAl09}, the VTSI analysis outputs such values as displacement, accelerations, and contact forces (as well as velocities) for all the model components. In practice, the obtained values have to be compared with the code provisions. 

\subsection{Realistic model}  \label{sec:moreNumResults:real}

\begin{figure}[H]
\centering
\setcounter{subfigure}{0}
\captionsetup[subfigure]{justification=centering}
\begin{tabular}{l@{\hskip 0.8cm}r}
\subfloat[Vertical bridge displacements] 
    {\label{fig:real:brDisplZ} 
    \includegraphics{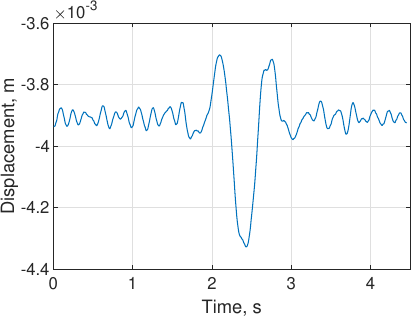}} &
\subfloat[Transverse bridge displacements] 
    {\label{fig:real:brDisplY} 
    \includegraphics{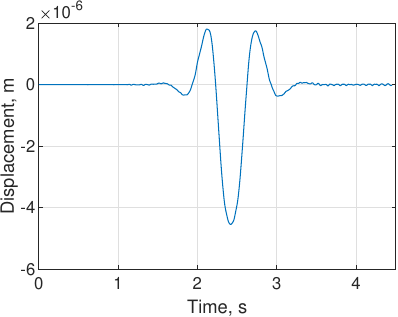}}
\end{tabular}
\caption{
    \textbf{\textit{Model:}} 
    Realistic vehicle passing over a 15-span bridge (Figure \ref{fig:modelPlanReal}) with speed $100\ \si{m/s}$. 
    \textbf{\textit{Algorithmic feature tested:}} 
    The Generalized-$\alpha$ scheme (Section \ref{sec:genAlpha}) with slight numerical damping ($\rhoinf = 0.9$).
    \textbf{\textit{Results:}} 
    Bridge displacements at the middle of the structure.
}
\label{fig:real:brDispl}
\end{figure}

\begin{figure}[H]
\centering
\setcounter{subfigure}{0}
\captionsetup[subfigure]{justification=centering}
\begin{tabular}{c}
\subfloat[Vertical bridge accelerations] 
    {\label{fig:real:brAccZ} 
    \includegraphics{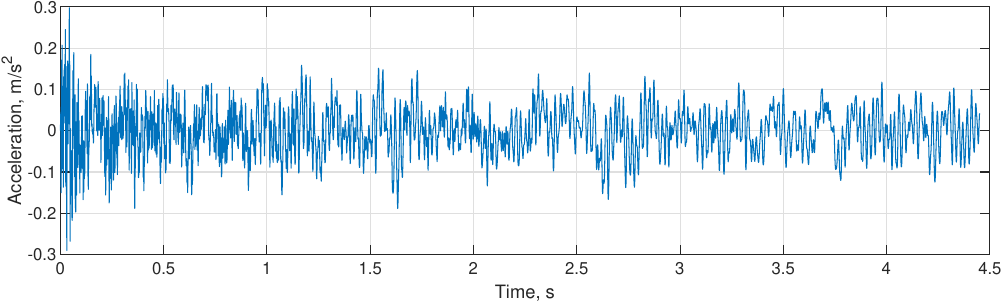}} \\
\subfloat[Transverse bridge accelerations] 
    {\label{fig:real:brAccY} 
    \includegraphics{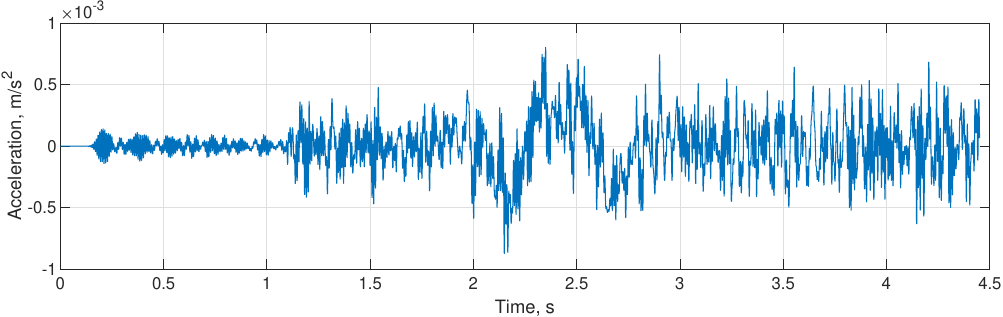}}
\end{tabular}
\caption{
    \textbf{\textit{Model:}} 
    Realistic vehicle passing over a 15-span bridge (Figure \ref{fig:modelPlanReal}) with speed $100\ \si{m/s}$. 
    \textbf{\textit{Algorithmic feature tested:}} 
    The Generalized-$\alpha$ scheme (Section \ref{sec:genAlpha}) with slight numerical damping ($\rhoinf = 0.9$).
    \textbf{\textit{Observations:}} 
    The frequency of bridge oscillations increases once realistic track and vehicle models are introduced (compare, for example, with the simplified model's results in Figure \ref{fig:NURBS_GenAl_corr:a}). However, the track provides additional damping to the system and redistributes forces, resulting in smaller amplitude changes. The use of the Generalized-$\alpha$ algorithm helps to mitigate spurious oscillations in this realistic model with many degrees of freedom. 
}
\label{fig:real:brAcc}
\end{figure}

\begin{figure}[H]
\centering
\setcounter{subfigure}{0}
\captionsetup[subfigure]{justification=centering}
\begin{tabular}{l@{\hskip 0.8cm}r}
\subfloat[Vertical car body displacements] 
    {\label{fig:real:carDisplZ} 
    \includegraphics{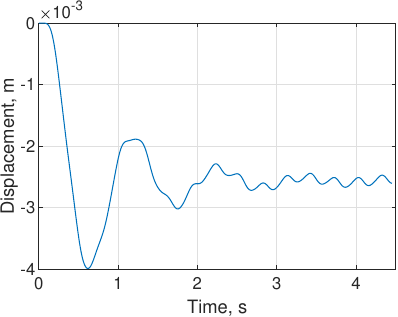}} &
\subfloat[Transverse car body displacements] 
    {\label{fig:real:carDisplY} 
    \includegraphics{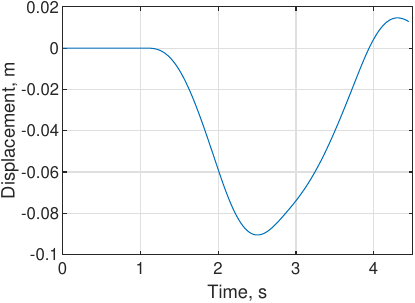}}
\end{tabular}
\caption{
    \textbf{\textit{Model:}} 
    Realistic vehicle passing over a 15-span bridge (Figure \ref{fig:modelPlanReal}) with speed $100\ \si{m/s}$. 
    \textbf{\textit{Algorithmic feature tested:}} 
    The Generalized-$\alpha$ scheme (Section \ref{sec:genAlpha}) with slight numerical damping ($\rhoinf = 0.9$).
    \textbf{\textit{Results:}} 
    Train car body displacements.
}
\label{fig:real:carDispl}
\end{figure}

\begin{figure}[H]
\centering
\setcounter{subfigure}{0}
\captionsetup[subfigure]{justification=centering}
\begin{tabular}{l@{\hskip 0.8cm}r}
\subfloat[Vertical car body accelerations] 
    {\label{fig:real:carAccZ} 
    \includegraphics{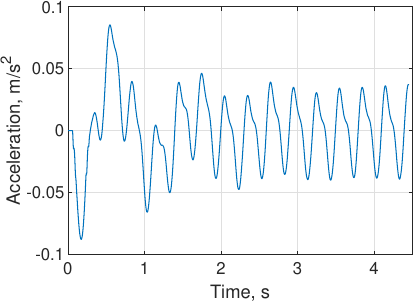}} &
\subfloat[Transverse car body accelerations] 
    {\label{fig:real:carAccY} 
    \includegraphics{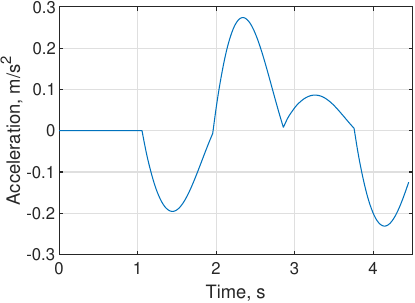}}
\end{tabular}
\caption{
    \textbf{\textit{Model:}} 
    Realistic vehicle passing over a 15-span bridge (Figure \ref{fig:modelPlanReal}) with speed $100\ \si{m/s}$. 
    \textbf{\textit{Algorithmic feature tested:}} 
    The Generalized-$\alpha$ scheme (Section \ref{sec:genAlpha}) with slight numerical damping ($\rhoinf = 0.9$).
    \textbf{\textit{Results:}} 
    Train car body accelerations.
}
\label{fig:real:carAcc}
\end{figure}

\begin{figure}[H]
\centering
\setcounter{subfigure}{0}
\captionsetup[subfigure]{justification=centering}
\begin{tabular}{c}
\subfloat[Vertical contact forces for the first wheelset] 
    {\label{fig:real:railTrainCF_Z} 
    \includegraphics{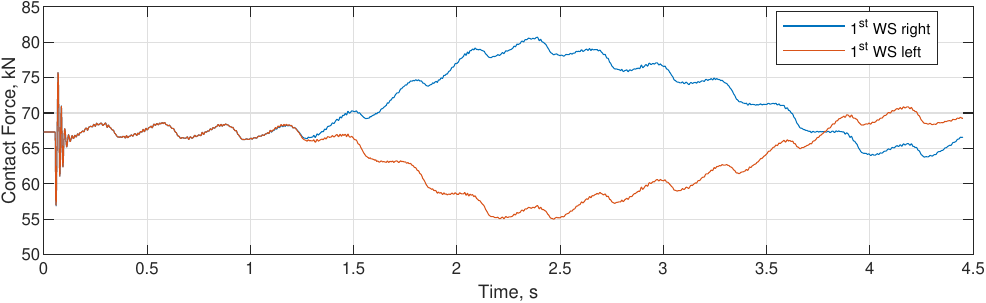}} \\
\subfloat[Vertical accelerations of a wheelset contact point (right wheel of the first wheelset)] 
    {\label{fig:real:CPacc_Z} 
    \includegraphics{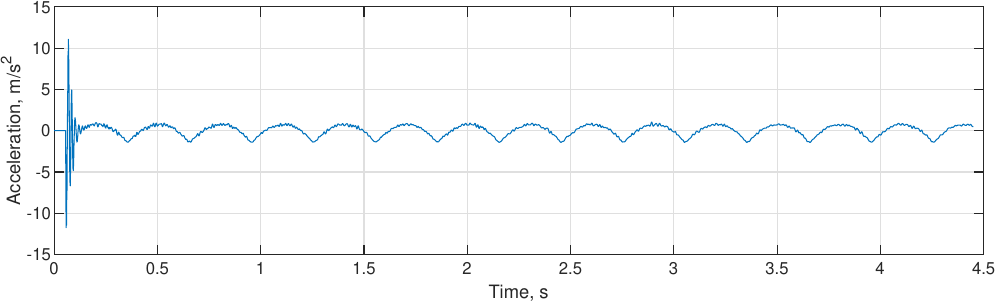}}
\end{tabular}
\caption{
    \textbf{\textit{Model:}} 
    Realistic vehicle passing over a 15-span bridge (Figure \ref{fig:modelPlanReal}) with speed $100\ \si{m/s}$. 
    \textbf{\textit{Algorithmic feature tested:}} 
    The Generalized-$\alpha$ scheme (Section \ref{sec:genAlpha}) with slight numerical damping ($\rhoinf = 0.9$).
    \textbf{\textit{Observations:}} 
    a) vertical contact force increases at the right wheel and decreases at the left wheel of the first wheelset while the train is passing over the curved spans; b) vertical accelerations of a wheel contact point are obtained. Acceleration discontinuity is observed when the vehicle enters the bridge, similar to the simplified model (Figure \ref{fig:NURBS_GenAl09:b}).
}
\label{fig:real:CPacc_force_Z}
\end{figure}

\begin{figure}[H]
\centering
\setcounter{subfigure}{0}
\captionsetup[subfigure]{justification=centering}
\begin{tabular}{c@{\hskip 0.8cm}c}
\subfloat[Transverse contact forces for all the wheelset] 
    {\label{fig:real:railTrainCF_Y:a} 
    \includegraphics{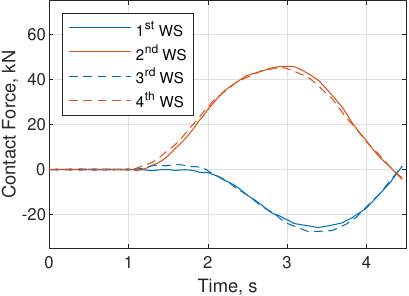}} &
\subfloat[A sum of transverse contact forces and comparison with a theoretical value] 
    {\label{fig:real:railTrainCF_Y_sum} 
    \includegraphics{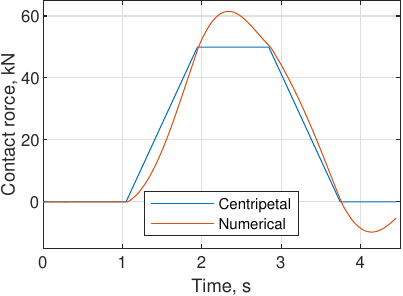}} \\
\multicolumn{2}{c} {
\subfloat[Transverse accelerations of a wheelset contact point (right wheel of the first wheelset)] 
    {\label{fig:real:CPacc_Y} \includegraphics{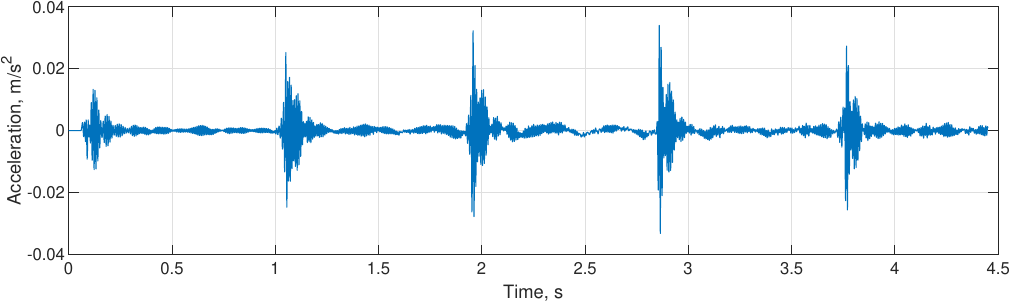}} }
\end{tabular}
\caption{
    \textbf{\textit{Model:}} 
    Realistic vehicle passing over a 15-span bridge (Figure \ref{fig:modelPlanReal}) with speed $100\ \si{m/s}$. 
    \textbf{\textit{Algorithmic feature tested:}} 
    The Generalized~-~$\alpha$ scheme (Section \ref{sec:genAlpha}) with slight numerical damping ($\rhoinf = 0.9$).
    \textbf{\textit{Observations:}} 
    Contact forces between wheels and rails are obtained.
    a) transverse contact forces for all the wheelsets; 
    b) sum of contact forces resulting from the centrifugal acceleration of the train compared with the theoretical value of the centripetal force (\textit{Note}: The difference between the numerical result and the theoretical is due to the fact that centripetal force is calculated for a single mass element only, not accounting for the effects of multiple bodies and suspensions);
    c) transverse accelerations of a wheel contact point.
}
\label{fig:real:railTrainCF_Y}
\end{figure}

\begin{figure}[H]
\centering
\setcounter{subfigure}{0}
\captionsetup[subfigure]{justification=centering}
\begin{tabular}{c}
\subfloat[Vertical contact force] 
    {\label{fig:real:railBrCF_Z} 
    \includegraphics{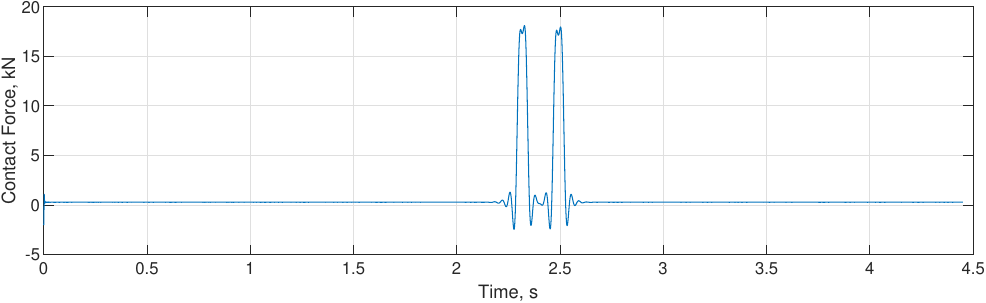}} \\
\subfloat[Transverse contact force] 
    {\label{fig:real:railBrCF_Y} 
    \includegraphics{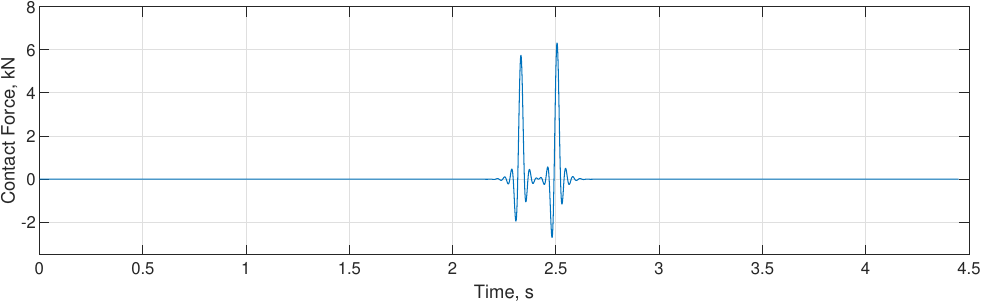}}
\end{tabular}
\caption{
    \textbf{\textit{Model:}} 
    Realistic vehicle passing over a 15-span bridge (Figure \ref{fig:modelPlanReal}) with speed $100\ \si{m/s}$. 
    \textbf{\textit{Algorithmic feature tested:}} 
    The Generalized-$\alpha$ scheme (Section \ref{sec:genAlpha}) with slight numerical damping ($\rhoinf = 0.9$).
    \textbf{\textit{Observations:}} 
    Contact forces between the right rail and the bridge are obtained in the middle of the structure. In the present study, the connection between the bridge structure and the rails is modeled using kinematic constraints. The use of The Generalized-$\alpha$ scheme and the rails formulated with NURBS allows for obtaining contact forces without spurious oscillations. 
}
\label{fig:real:railBrCF}
\end{figure}

A realistic train model from Section \ref{sec:trainReal} and realistic bridge and track models from Section \ref{sec:bridgeReal} are used in this section. A one-car train (Figure \ref{fig:trainModelReal}) is passing over a 15-span bridge (Figure \ref{fig:modelPlanReal}) with speed $100\ \si{m/s}$, reaching the end of the bridge at time $4.5\ \si{s}$.  The results in Figures \ref{fig:real:brDispl} - \ref{fig:real:railBrCF} are obtained using the Generalized-$\alpha$ scheme with numerical dissipation. It was observed that for models with many degrees of freedom, some numerical damping is needed, even if other strategies, such as NURBS-based spatial discretization and temporal constraint projection, are employed, since, in such models, the oscillations arise not only due to the presence of the constraints but also due to other high-frequency components in the solution. Therefore, if such a method as constraint projection is used, a very small time step is required to avoid the oscillations, which is not practical in realistic applications.

Eurocode EN 1990 (2002)\cite{EN1990} limits the value of vertical bridge deck acceleration to $3.5\ \si{m/s^2}$ (for a ballasted track). The values presented in Figure \ref{fig:real:brAccZ} are well below the limit. The low value can be attributed to the fact that a train composed of only a single car was considered in this example to illustrate the modeling approaches described in previous sections. A more significant response, including a resonance, is expected when the bridge is analyzed under the trains with a standard number of cars. To provide a "very good" level of comfort to passengers, EN 1990\cite{EN1990} limits vertical accelerations of a train car to $1\ \si{m/s^2}$. The values in Figure \ref{fig:real:carAccZ} are below this limit.

\section{Concluding remarks}    \label{sec:conclusions}


An algorithm for three-dimensional dynamic vehicle-track-structure interaction (VTSI) analysis has been discussed in this paper. It builds on our previous work on a two-dimensional VTSI algorithm \cite{Fedorova2017} and includes three main new features: 1) a corotational approach to modeling vehicle dynamics, 2) NURBS-based discretization of the vehicle path, 3) two alternative approaches to solving an index-3 DAE system that are possible due to higher order continuity provided by NURBS. 
 
Employing the corotational approach, the train dynamics are conveniently described with respect to Frenet frames along the train's path on the bridge. These Frenet frames serve as corotational frames. While displacements and rotations of the train can be large, deformations within these corrotational frames are assumed to be small, which allows the use of a linear (time-varying) dynamic model. First, a simplified VTSI model is discussed in this paper for clarity of demonstration, followed by an example of a realistic model. Complete multibody train framework, as well as interacting bridge and track models, are discussed in Reference \cite{FedorovaThesis2023}.

While solving VTSI equations as an index-3 DAE system and employing standard finite elements to model the track, spurious oscillations occur in obtained contact forces and accelerations without appropriate algorithmic measures. Moreover, if straight finite elements are used to approximate a curved bridge, spurious oscillations are amplified in the transverse direction due to curvature and rotation discontinuities. We employ isogeometric analysis to model the bridge using NURBS. This helps mitigate oscillations in the transverse direction. However, NURBS shape functions are twice continuously differentiable only inside the bridge. This results in a discontinuity in wheel accelerations at the beginning and end of the bridge. Correction steps are proposed to avoid this discontinuity. Overall, three alternative approaches to solving the VTSI system are discussed (see Figure \ref{fig:summary} for the summary).

The higher order continuity due to the use of NURBS allows for two alternative approaches for solving the VTSI system.
First, if the VTSI system is reduced to index-1 DAE (when constraints are enforced at acceleration level), only initial wheel velocity and acceleration need to be corrected to avoid the discontinuity at the beginning of the bridge and the Newmark scheme can be used for time integration (while avoiding numerical damping). This correction ensures consistent initial conditions at all three levels (displacement, velocity, and acceleration) of constraints. When displacement constraints are not enforced explicitly, drift in displacement compatibility develops. The observed drift is not significant, but a periodic wheel displacement correction based on the kinematic constraint can be used as needed. Second, it has also been demonstrated that index-3 VTSI equations can be solved using the Newmark scheme if wheel velocities and accelerations are corrected at each time step using a simple non-iterative constraint projection method. 

The proposed approach allows for a convenient description of vehicle dynamics through corotational frames and for solving three-dimensional VTSI as an index-3 DAE system with or without numerical damping. The algorithm can be applied to other multibody systems (with large displacements and rotations, but small deformations) coupled to deformable structures by time-varying kinematic constraints. The use of kinematic constraints leads to the modularity of the proposed VTSI algorithm and its straightforward integration into existing finite element software.

\FloatBarrier


\section*{Acknowledgments}
The first author gratefully acknowledges a Ph.D. scholarship from LARSA Inc., Melville, NY.


\appendix
\section{Simplified vehicle matrices} \label{sec:appendix:train}

\setcounter{MaxMatrixCols}{20}
\arraycolsep = 1pt       

This appendix includes matrices for a simplified vehicle model presented in Section \ref{sec:trainSimple}. Here, $\Mt$, $\Ct(t)$, and $\Kt(t)$ are the mass, damping, and stiffness matrices. $\Ltr$ is the kinematic constraint matrix used to constrain wheels to the rails. $\lambdaB$ is a vector of constraint forces in transverse, vertical, and rolling directions. $\Pt$ is the load vector. 

\begin{equation} \label{train3D_members}
\Mt = \left[ {\begin{array}{*{20}{c}}
{\mw+\mc}&{}&{{-\mc}{l_0}}&{}\\
{}&{{\mw}}&{}&{}\\
{{-\mc}{l_0}}&{}&{{\mc}{l_0}^2 + {\Iw} + {\Ic}}&{}\\
{}&{}&{}&{{\mc}}
\end{array}} \right]
\end{equation}

\begin{equation}
\Ct = 
2 \mw {\Tw}^\T \omegaFhat \Tw + 2 \mc {\Tc}^\T \omegaFhat \Tc = 
\left[ {\begin{array}{*{20}{c}}
{}&{ - 2{\mw}\omegaF_1}&{}&{ - 2{\mc}\omegaF_1}\\
{2{\mw}\omegaF_1}&{}&{}&{}\\
{}&{}&{}&{  2{\mc}\omegaF_1 {l_0}}\\
{2{\mc}\omegaF_1}&{}&{-2{\mc}\omegaF_1{l_0}}&{}
\end{array}} \right]
\end{equation}

\begin{equation}
\begin{gathered}
\Kt = 
\mw {\Tw}^\T \omegaFhatdot \Tw + \mc {\Tc}^\T \omegaFhatdot \Tc +
\mw {\Tw}^\T {\omegaFhat}^2 \Tw + \mc {\Tc}^\T {\omegaFhat}^2 \Tc= 
\\
\left[ {\begin{array}{*{20}{c}}
{ - \left( {{\mw} + {\mc}} \right)\left( {\omegaF_1}^2 + {\omegaF_3}^2 \right)}
&{ - {\mw}\left( \omegaFdot_1 - \omegaF_2 \omegaF_3 \right)}
&{   {\mc}\left( {\omegaF_1}^2 + {\omegaF_3}^2 \right){l_0}}
&{ - {\mc}\left( \omegaFdot_1 - \omegaF_2 \omegaF_3 \right)}
\\
{{\mw}\left( {\omegaFdot_1 + \omegaF_2 \omegaF_3} \right)}
&{ - {\mw}\left( {\omegaF_1}^2 + {\omegaF_2}^2 \right) + {\ks}}
&0
&{ - {\ks}}
\\
{   {\mc}\left( {\omegaF_1}^2 + {\omegaF_3}^2 \right){l_0}}
&0
&{ - {\mc}\left( {\omegaF_1}^2 + {\omegaF_3}^2 \right){l_0}^2}
&{   {\mc}\left( {\omegaFdot_1 - \omegaF_2 \omegaF_3} \right){l_0}}
\\
{{\mc}\left( {\omegaFdot_1 + \omegaF_2 \omegaF_3 } \right)}
&{ - {\ks}}
&{ - {\mc}\left( {\omegaFdot_1 + \omegaF_2 \omegaF_3} \right){l_0}}
&{ - {\mc}\left( {\omegaF_1}^2 + {\omegaF_2}^2 \right) + {\ks}}
\end{array}} \right]
\end{gathered}
\end{equation}

\begin{equation}
\Ltr = \left[ {\begin{array}{*{20}{c}}
{ - 1}&0&0\\
0&{ - 1}&0\\
0&0&{ - 1}\\
0&0&0
\end{array}} \right]
\end{equation}

\begin{equation}
\lambdaB = \left( {\begin{array}{*{20}{c}}
{{\lambda_y}}\\
{{\lambda_z}}\\
{{\lambda_{{\theta _x}}}}
\end{array}} \right) 
\end{equation}

\begin{equation}    \label{eq:Pt}
\begin{gathered}
\Pt = 
- \left( \mw {\Tw}^\T + \mc {\Tc}^\T \right) {\RF}^\T  \xFddot 
- \mc {\Tc}^\T \omegaFhatdot 
    \left( {\begin{array}{*{20}{c}}
    0\\
    0\\
    l_0
    \end{array}} \right)
- \mc {\Tc}^\T {\omegaFhat}^2
    \left( {\begin{array}{*{20}{c}}
    0\\
    0\\
    l_0
    \end{array}} \right)
\\
- (\mw g {\Tw}^\T+ \mc g {\Tc}^\T) ({\RF(t)}^\T - {\RF(0)}^\T) 
\left( {\begin{array}{*{20}{c}}
    0\\
    0\\
    1
    \end{array}} \right)
\end{gathered}
\end{equation}

The vehicle matrices listed above are used in the simplified vehicle equations of motion \eqref{eq:t:eqofmotion}. Complete derivations of the equations of motion can be found in Fedorova \cite{FedorovaThesis2017}.





\bibliographystyle{elsarticle-num-names} 
\bibliography{cas-refs}





\end{document}